\documentclass[11pt]{article}
\usepackage{graphicx}
\usepackage{amsbsy}
\usepackage{amstext}
\usepackage{amssymb}
\usepackage{amsthm}
\usepackage{amsfonts}
\usepackage{amsmath}
\usepackage{subfig}
\usepackage{makeidx}
\usepackage{epstopdf}
\usepackage{dcolumn}
\usepackage{booktabs}
\usepackage[ruled,vlined]{algorithm2e}
\usepackage{xcolor}
\usepackage{times}
\usepackage{fullpage}
\usepackage[square,sort,comma,numbers]{natbib}

\newcommand{\svskip}{\vspace{1.75mm}}

\newcommand{\tr}{\mathop{\rm tr}\nolimits}

\newcommand{\conv}{\mathop{\rm conv}\nolimits}
\newcommand{\dom}{\mathop{\rm dom}\nolimits}

\newcommand{\dist}{\mathop{\rm dist}\nolimits}
\newcommand{\prox}{\mathop{\rm prox}\nolimits}
\newcommand{\argmin}{\mathop{\rm argmin}\nolimits}
\newcommand{\amp}{\mathop{\;\:}\nolimits}

\newcommand{\ba}{\boldsymbol{a}}
\newcommand{\bb}{\boldsymbol{b}}
\newcommand{\bc}{\boldsymbol{c}}

\newcommand{\be}{\boldsymbol{e}}

\newcommand{\bp}{\boldsymbol{p}}

\newcommand{\bs}{\boldsymbol{s}}

\newcommand{\bu}{\boldsymbol{u}}
\newcommand{\bv}{\boldsymbol{v}}
\newcommand{\bw}{\boldsymbol{w}}
\newcommand{\bx}{\boldsymbol{x}}
\newcommand{\by}{\boldsymbol{y}}
\newcommand{\bz}{\boldsymbol{z}}
\newcommand{\bA}{\boldsymbol{A}}

\newcommand{\bD}{\boldsymbol{D}}
\newcommand{\bE}{\boldsymbol{E}}

\newcommand{\bG}{\boldsymbol{G}}

\newcommand{\bI}{\boldsymbol{I}}

\newcommand{\bL}{\boldsymbol{L}}
\newcommand{\bM}{\boldsymbol{M}}
\newcommand{\bN}{\boldsymbol{N}}

\newcommand{\bQ}{\boldsymbol{Q}}

\newcommand{\bS}{\boldsymbol{S}}

\newcommand{\bU}{\boldsymbol{U}}
\newcommand{\bV}{\boldsymbol{V}}
\newcommand{\bW}{\boldsymbol{W}}
\newcommand{\bX}{\boldsymbol{X}}
\newcommand{\bY}{\boldsymbol{Y}}
\newcommand{\bZ}{\boldsymbol{Z}}

\newcommand{\blambda}{\boldsymbol{\lambda}}

\newcommand{\bSigma}{\boldsymbol{\Sigma}}

\newcommand{\bzero}{\boldsymbol{0}}

\newtheorem{proposition}{Proposition}

\SetAlCapSty{xAlCapSty}

\SetCommentSty{xCommentSty}
\SetKwComment{Comment}{\# }{}

\SetNlSty{mynlfont}{}{}

\SetSideCommentRight

\DontPrintSemicolon

\RestyleAlgo{algoruled}

\begin{document}

\title{Proximal Distance Algorithms: Theory and Examples}
\author{Kevin L. Keys$^{1,2}$ \\
Hua Zhou$^{3}$ \\
Kenneth Lange$^{2,4,5}$ \\
\\
Department of Medicine$^{1}$ \\
University of California, San Francisco, CA 94158 \\
Departments of Biomathematics$^{2}$, Biostatistics$^{3}$, Human Genetics$^{4}$, and Statistics$^{5}$ \\
University of California, \\
Los Angeles, CA 90095-1766}

\maketitle

\begin{abstract}
Proximal distance algorithms combine the classical penalty method of constrained minimization with distance majorization. If $f(\bx)$ is the loss function, and $C$ is the constraint set in a constrained minimization problem, then the proximal distance principle mandates minimizing the penalized loss $f(\bx)+\frac{\rho}{2}\dist(\bx,C)^2$ and following the solution $\bx_{\rho}$ to its limit as $\rho$ tends to $\infty$. At each iteration the squared Euclidean distance $\dist(\bx,C)^2$ is majorized by the spherical quadratic $\|\bx-P_C(\bx_k)\|^2$, where $P_C(\bx_k)$ denotes the projection of the current iterate $\bx_k$ onto $C$. The minimum of the surrogate function $f(\bx)+\frac{\rho}{2}\|\bx-P_C(\bx_k)\|^2$ is given by the proximal map $\prox_{\rho^{-1}f}[P_C(\bx_k)]$. The next iterate $\bx_{k+1}$ automatically decreases the original penalized loss for fixed $\rho$. Since many explicit projections and proximal maps are known, it is straightforward to derive and implement novel optimization algorithms in this setting. These algorithms can take hundreds if not thousands of iterations to converge, but the simple nature of each iteration makes proximal distance algorithms competitive with traditional algorithms. For convex problems, proximal distance algorithms reduce to proximal gradient algorithms and therefore enjoy well understood
convergence properties. For nonconvex problems, one can attack convergence by invoking
Zangwill's theorem. Our numerical examples demonstrate the utility of proximal distance algorithms in various high-dimensional settings, including a) linear programming, b) constrained least squares, c) projection to the closest kinship matrix, d) projection onto a second-order cone constraint, e) calculation of Horn's copositive matrix index, f) linear complementarity programming, and g) sparse principal components analysis. The proximal distance algorithm in each case is competitive or
    superior in speed to traditional methods such as the interior point method and the alternating direction method of multipliers (ADMM). \\
\svskip 
\textit{Math Subject Classifications:} 90C59, 90C26, 65K05
\end{abstract}

\section{Introduction}
The solution of constrained optimization problems is part science and part art. As mathematical scientists explore the largely uncharted territory of high-dimensional nonconvex problems, it is imperative to consider new methods. The current paper studies a class of optimization algorithms that combine Courant's penalty method of optimization \citep{beltrami1970algorithmic,courant1943variational} 
with the notion of a proximal operator \citep{bauschke2011convex,moreau1962fonctions,parikhboyd2013}. 
The classical penalty method turns constrained minimization of a function $f(\bx)$ over a closed set $C$ into unconstrained minimization. The general idea is to seek the minimum point of a penalized version $f(\bx)+\rho q(\bx)$ of $f(\bx)$, where the penalty $q(\bx)$ is nonnegative and vanishes precisely on $C$. If one follows the solution vector $\bx_\rho$ as $\rho$ tends to $\infty$, then in the limit one recovers the constrained solution. The penalties of choice in the current paper are squared Euclidean distances $\dist(\bx,C)^2=\inf_{\by \in C}\|\bx-\by\|^2$.

The formula
\begin{eqnarray}
\prox_f(\by) & = & \argmin_{\bx}\left[f(\bx)+\frac{1}{2}\|\bx-\by\|^2\right]
\label{prox_operator}
\end{eqnarray} 
defines the proximal map of a function $f(\bx)$. Here $\|\cdot\|$ is again the
standard Euclidean norm, and $f(\bx)$ is typically assumed to be closed
and convex. Projection onto a closed convex set $C$ is realized by choosing $f(\bx)$ to be the $0/\infty$ indicator $\delta_C(\bx)$ of $C$. It is possible to drop the convexity assumption if $f(\bx)$ is nonnegative or coercive. In so doing, $\prox_f(\by)$ may become multi-valued. For example, the minimum distance from a nonconvex set to an exterior point may be attained at multiple boundary points. The point 
$\bx$ in the definition (\ref{prox_operator}) can be restricted to a subset $S$ of Euclidean space by replacing $f(\bx)$ by $f(\bx) + \delta_S(\bx)$, where $\delta_S(\bx)$ is the indicator of $S$.
 
One of the virtues of exploiting proximal operators is that they have
been thoroughly investigated. For a large number of functions $f(\bx)$,
the map $\prox_{cf}(\by)$ for $c > 0$ is either given by an exact formula or 
calculable by an efficient algorithm. The known formulas tend to be highly accurate. This is a plus because the classical penalty method suffers from ill conditioning for large values of the penalty 
constant. Although the penalty method seldom delivers exquisitely accurate solutions, moderate accuracy suffices for many problems. 

There are ample precedents in the optimization literature for the 
proximal distance principle. Proximal gradient algorithms have been employed for many years in many contexts, including projected Landweber, alternating projection onto the intersection of two or more closed convex sets, the alternating-direction method of multipliers (ADMM), and fast iterative shrinkage thresholding algorithms (FISTA) 
\citep{beck2009fast,combettes2011proximal,landweber1951iteration}.
Applications of distance majorization are more recent \citep{chi2014distance,lange2014proximal,xu2017generalized}. The overall strategy consists of replacing the distance penalty $\dist(\bx,C)^2$
by the spherical quadratic $\|\bx-\by_k\|^2$, where $\by_k$ is the projection of the $k$th iterate $\bx_k$ onto $C$. To form the next iterate, one then sets 
\begin{eqnarray*}
\bx_{k+1}& = & \prox_{\rho^{-1}f}(\by_k)  \quad \text{with} \quad
\by_k \amp = \amp P_C(\bx_k).
\end{eqnarray*}
The MM (majorization-minimization) principle guarantees that $\bx_{k+1}$
decreases the penalized loss. We call the combination of Courant's penalty method with distance majorization the \emph{proximal distance principle}.
Algorithms constructed according to the principle are \emph{proximal distance algorithms}. 

The current paper extends and deepens our previous preliminary treatments of the proximal distance principle. Details of implementation such as Nesterov acceleration matter in performance. We have found that squared distance penalties tend to work better than exact penalties. In the presence of convexity, it is now clear that every proximal distance algorithm reduces to a proximal gradient algorithm. Hence, convergence analysis can appeal to a venerable body of convex theory. This does not imply that the proximal distance algorithm is limited to convex problems. In fact, its most important applications may well be to nonconvex problems. A major focus of this paper is on practical exploration of the proximal distance algorithm. 

The current paper also presents some fresh ideas. Among the innovations are: a) recasting proximal distance algorithms with convex losses as concave-convex programs, b) providing a theoretical convergence analysis for nonconvex proximal distance algorithms, c) demonstrating the virtue of folding constraints into the domain of the loss, and d) treating in detail seven interesting examples. It is noteworthy that our new convergence 
theory covers more general MM algorithms. It is our sincere hope to enlist other mathematical scientists in expanding and clarifying this promising line of research. A clear exposition of the known facts seems like the logical place to start.

To that end, we now outline the remainder of our paper. Section \ref{sec:derivation} briefly sketches the underlying MM principle.
We then show how to construct proximal distance algorithms from the MM principle and distance majorization. The Section concludes with the
derivation a few broad categories proximal distance algorithms. Section \ref{sec:conv_accel} covers convergence theory for convex problems, while Section \ref{sec:nonconvex_convergence} provides a more general treatment of convergence for nonconvex problems. Section \ref{sec:examples} discusses our numerical experiments on various convex and nonconvex problems. We close by indicating some future research directions in Section \ref{sec:discussion}.

\section{Derivation}
\label{sec:derivation}

The derivation of our proximal distance algorithms exploits the majorization-minimiza\-tion (MM) principle \citep{hunter04,lange10}.
In minimizing a function $f(\bx)$, the MM principle exploits a surrogate function $g(\bx \mid \bx_k)$ that majorizes $f(\bx)$ around the current iterate $\bx_k$. Majorization mandates both domination $g(\bx \mid \bx_k) \ge f(\bx)$ for all feasible $\bx$ and tangency $g(\bx_k \mid \bx_k) = f(\bx_k)$ at the anchor $\bx_{k}$.  If $\bx_{k+1}$  minimizes $g(\bx \mid \bx_k)$, then the descent property $f(\bx_{k+1}) \le f(\bx_k)$ follows from the string of inequalities and equalities
\begin{eqnarray*}
f(\bx_{k+1}) &\le & g(\bx_{k+1} \mid \bx_k) \amp \le \amp g(\bx_k \mid \bx_k)
\amp = \amp  f(\bx_k) .
\end{eqnarray*}
Clever selection of the surrogate $g(\bx \mid \bx_{k+1})$ can lead to a simple algorithm with an explicit update that requires little computation per iterate. The number of iterations until convergence of an MM algorithm depends on how tightly $g(\bx \mid \bx_k)$ hugs
$f(\bx)$. Constraint satisfaction is built into any MM algorithm. 
If maximization of $f(\bx)$ is desired, then the objective $f(\bx)$ should dominate the surrogate $g(\bx \mid \bx_k)$ subject to the tangency condition. The next iterate $\bx_{k+1}$ is then chosen
to maximize $g(\bx \mid \bx_k)$. The minorization-maximization version of the MM principle guarantees the ascent property.

The constraint set $C$ over which the loss $f(\bx)$ is minimized
can usually be expressed as an intersection $\cap_{i=1}^m C_i$
of closed sets. It is natural to define the penalty
\begin{eqnarray*}
q(\bx) & = &  \frac{1}{2} \sum_{i=1}^m \alpha_i \dist(\bx,C_i)^2
\end{eqnarray*}
using a convex combination of the squared distances. The
neutral choice $\alpha_i=\frac{1}{m}$ is one we prefer in practice.
Distance majorization gives the surrogate function
\begin{eqnarray*}
g_\rho(\bx \mid \bx_k) & = & f(\bx)+\frac{\rho}{2} \sum_{i=1}^m \alpha_i \|\bx-P_{C_i}(\bx_k)\|^2 \\
& = & f(\bx)+\frac{\rho}{2}\Big\|\bx-\sum_{i=1}^m \alpha_i P_{C_i}(\bx_k)\Big\|^2+c_k
\end{eqnarray*}
for an irrelevant constant $c_k$. If we put $\by_k=\sum_{i=1}^m \alpha_i P_{C_i}(\bx_k)$,
then by definition the minimum of the surrogate $g_\rho(\bx \mid \bx_k)$ occurs at the proximal point
\begin{eqnarray}
\bx_{k+1} & = & \prox_{\rho^{-1}f}(\by_k). \label{algorithm_map}
\end{eqnarray} 
We call this MM algorithm the proximal distance algorithm.
The penalty $q(\bx)$ is generally smooth because 
\begin{eqnarray*}
    \nabla \frac{1}{2}\dist(\bx,C)^2 & = &\bx - P_C(\bx)
\end{eqnarray*}
at any point $\bx$ where the projection $P_C(\bx)$ is single valued \citep{borwein06,lange2016MM}. This is always true for convex sets
and almost always true for nonconvex sets. For the moment, we will
ignore the possibility that $P_C(\bx)$ is multi-valued.

For the special case of projection of an external point $\bz$ onto the intersection $C$
of the closed sets $C_i$, one should take $f(\bx)=\frac{1}{2}\|\bz-\bx\|^2$. The proximal distance iterates then obey the explicit formula
\begin{eqnarray*}
\bx_{k+1} & = & \frac{1}{1+\rho}(\bz+\rho \by_k).
\end{eqnarray*}
Linear programming with arbitrary convex constraints is another example. Here 
the loss is $f(\bx) = \bv^t\bx$, and the update reduces to
\begin{eqnarray*}
\bx_{k+1} & = & \by_k-\frac{1}{\rho}\bv.
\end{eqnarray*}
If the proximal map is impossible to calculate, but $\nabla f(\bx)$ is 
known to be Lipschitz with constant $L$, then one can substitute the standard majorization 
\begin{eqnarray*}
f(\bx) & \le & f(\bx_k)+\nabla f(\bx_k)^t(\bx-\bx_k)+\frac{L}{2}\|\bx-\bx_k\|^2
\end{eqnarray*}
for $f(\bx)$. Minimizing the sum of the loss majorization plus the penalty majorization
leads to the MM update
\begin{eqnarray}
\bx_{k+1} & = & \frac{1}{L+\rho}[-\nabla f(\bx_k)+L\bx_k+\rho \by_k] \nonumber \\
& = & \bx_k - \frac{1}{L+\rho}[\nabla f(\bx_k)+\rho \nabla q(\bx_k)]  .
\label{double_majorization_solution}
\end{eqnarray}
This is a gradient descent algorithm without an intervening proximal map.

In moderate-dimensional problems, local quadratic approximation of $f(\bx)$ can lead to a viable algorithm. For instance, in generalized linear statistical models, \cite{xu2017generalized} suggest replacing the observed information matrix by the expected information
matrix. The latter matrix has the advantage of being positive semidefinite. In our notation, if $\bA_k \approx d^2f(\bx_k)$, then an approximate quadratic surrogate is
\begin{eqnarray*}
f(\bx_k)+\nabla f(\bx_k)^t (\bx-\bx_k)+\frac{1}{2}(\bx-\bx_k)^t\bA_k(\bx-\bx_k)+\frac{\rho}{2}\|\bx-\by_k\|^2.
\end{eqnarray*}
The natural impulse is to update $\bx$ by the Newton step
\begin{eqnarray}
\bx_{k+1} & = & \bx_k-(\bA_k+\rho \bI)^{-1}[\nabla f(\bx_k)
-\rho \by_k]. \label{newton_update}
\end{eqnarray}
This choice does not necessarily decrease $f(\bx)$. Step halving or another form of backtracking restores the descent property. 

A more valid concern is the effort expended in matrix inversion. If $\bA_k$ is dense and constant, then extracting the spectral decomposition $\bV \bD \bV^t$ of $\bA$ reduces formula (\ref{newton_update}) to
\begin{eqnarray*}
\bx_{k+1} & = & \bx_k-\bV(\bD+\rho \bI)^{-1}\bV^t[\nabla f(\bx_k)
-\rho \by_k],
\end{eqnarray*}
which can be implemented as a sequence of matrix-vector multiplications. Alternatively, one can take just a few terms of the series
\begin{eqnarray*}
(\bA_k+\rho \bI)^{-1} & = & \rho^{-1}\sum_{j=0}^\infty 
(-\rho^{-1}\bA_k)^j
\end{eqnarray*}
when $\rho$ is sufficiently large. For a generalized linear model, parameter updating involves solving the linear system
\begin{eqnarray}
(\bZ^t\bW_{\!k}\bZ+\rho \bI)\bx & = & \bZ^t \bW_{\!k}^{1/2} \bv_k + \rho \by_k \label{least_squares_update}
\end{eqnarray}
for $\bW_{\!k}$ a diagonal matrix with positive diagonal entries. This
task is equivalent to minimizing the least squares criterion
\begin{eqnarray}
\left\|\begin{pmatrix} \bW_{\!k}^{1/2}\bZ \\
\sqrt{\rho}\bI \end{pmatrix}\bx -\begin{pmatrix} \bv_k
\\ \sqrt{\rho}\by_k \end{pmatrix}
\right\|^2. \label{johann_credit}
\end{eqnarray}
In the unweighted case, extracting the singular value decomposition $\bZ = \bU \bS \bV^T$ facilitates solving the system 
of equations (\ref{least_squares_update}). The svd decomposition is 
especially cheap if there is a substantial mismatch between the 
number rows and columns of $\bZ$. For sparse $\bZ$, the conjugate gradient algorithm adapted to least squares \citep{paige1982algorithm} is subject to much less ill conditioning than the standard conjugate gradient algorithm. Indeed, the algorithm LSQR and its sparse version LSMR \citep{fong2011lsmr} perform well even when the matrix $(\bZ^t \bW_{\!k}^{1/2}, \sqrt{\rho}\bI)^t$ is ill conditioned. 

The proximal distance principle also applies to unconstrained problems. 
For example, consider the problem of minimizing a
penalized loss $\ell(\bx)+ p(\bA \bx)$. The presence of the linear transformation
$\bA \bx$ in the penalty complicates optimization. The strategy of 
parameter splitting introduces a new variable $\by$ and minimizes 
$\ell(\bx)+ p(\by)$ subject to the constraint $\by= \bA \bx$. If $P_M(\bz)$ 
denotes projection onto the manifold 
\begin{eqnarray*}
M = \{\bz = (\bx,\by): \bA\bx = \by\},
\end{eqnarray*}
then the constrained problem can be solved approximately by minimizing the
function 
\begin{eqnarray*}
\ell(\bx)+p(\by)+\frac{\rho}{2}\dist(\bz,M)^2
\end{eqnarray*}
for large $\rho$. If $P_M(\bz_k)$ consists of  two subvectors $\bu_k$ and
$\bv_k$ corresponding to $\bx_k$ and $\by_k$, then the proximal distance
updates are
\begin{eqnarray*}
\bx_{k+1} & = & \prox_{\rho^{-1} \ell}(\bu_k) \quad \text{and} \quad
\by_{k+1} \amp = \amp \prox_{\rho^{-1} p}(\bv_k) .
\end{eqnarray*}

Given the matrix $\bA$ is $n \times p$, one can attack the projection by minimizing the function
\begin{eqnarray*}
q(\bx) & = & \frac{1}{2}\|\bx-\bu\|^2+\frac{1}{2}\|\bA \bx-\bv\|^2.
\end{eqnarray*}
This leads to the solution
\begin{eqnarray*}
\bx & = & (\bI_p + \bA^t \bA)^{-1}(\bA^t\bv+\bu) \quad \text{and}
\quad \by \amp = \amp \bA\bx.
\end{eqnarray*}
If $n < p$, then the Woodbury formula
\begin{eqnarray*}
(\bI_p+\bA^t \bA)^{-1} & =  &
\bI_p-\bA^t(\bI_n+\bA\bA^t)^{-1}\bA
\end{eqnarray*}
reduces the expense of matrix inversion.

Traditionally, convex constraints have been posed as inequalities 
$C = \{\bx: a(\bx) \le t\}$. \cite{parikhboyd2013} point out how to project onto such sets. The relevant Lagrangian for projecting an external point $\by$ amounts to
\begin{eqnarray*}
\mathcal{L}(\bx,\lambda) & = & \frac{1}{2}\|\by-\bx\|^2
+ \lambda [a(\bx)-t] 
\end{eqnarray*}
with $\lambda \ge 0$. The corresponding stationarity condition
\begin{eqnarray}
{\bf 0} & = & \bx-\by + \lambda \nabla a(\bx), \label{boyd_bisect}
\end{eqnarray}
can be interpreted as $a[\prox_{\lambda a}(\by)] = t$. One 
can solve this one-dimensional equation for $\lambda$ by bisection. Once $\lambda$ is available, $\bx = \prox_{\lambda a}(\by)$ is available as well. \cite{parikhboyd2013} note that the value $a[\prox_{\lambda a}(\by)]$ is decreasing in $\lambda$. One can verify their claim by implicit differentiation of equation (\ref{boyd_bisect}). This gives
\begin{eqnarray*}
\frac{d}{d \lambda}\bx & = & - [\bI+\lambda d^2 a(\bx)]^{-1} \nabla a(\bx)
\end{eqnarray*}
and consequently the chain rule inequality
\begin{eqnarray*}
\frac{d}{d \lambda} a[\prox_{\lambda a}(\by)] & = & - da(\bx)[\bI+\lambda d^2 a(\bx)]^{-1} \nabla a(\bx) \amp \le \amp 0.
\end{eqnarray*}

\section{Convergence: Convex Case}
\label{sec:conv_accel}

In the presence of convexity, the proximal distance algorithm reduces to a proximal gradient algorithm.  This follows from the representation
\begin{eqnarray*}
\by & = & \sum_{i=1}^m \alpha_i P_{C_i}(\bx) \amp = \amp \bx 
- \sum_{i=1}^m \alpha_i \Big[\bx-P_{C_i}(\bx)\Big] \amp = \amp \bx - \nabla q(\bx)
\end{eqnarray*}
involving the penalty $q(\bx)$. Thus, the proximal distance algorithm can 
be expressed as 
\begin{eqnarray*}
\bx_{k+1} & = & \prox_{\rho^{-1}f}[\bx_k - \nabla q(\bx_k)].
\end{eqnarray*}
In this regard, there is the implicit assumption that $\nabla q(\bx)$ is Lipschitz
with constant 1. This is indeed the case. According to the Moreau decomposition
\citep{bauschke2011convex}, for a single closed convex set $C$
\begin{eqnarray*}
\nabla q(\bx) & = & \bx - P_C(\bx) \amp = \amp \prox_{\delta_C^\star}(\bx) ,
\end{eqnarray*}
where $\delta_C^\star(\bx)$ is the Fenchel conjugate of the indicator function 
\begin{eqnarray*}
\delta_C(\bx) & = & \begin{cases} 0 & \bx \in C \\ \infty & \bx \not\in C .\end{cases}
\end{eqnarray*}
Because proximal operators of closed convex functions are nonexpansive \citep{bauschke2011convex}, 
the result follows for a single set. For the 
general penalty $q(\bx)$ with $m$ sets, the Lipschitz constants are scaled by the
convex coefficients $\alpha_i$ and added to produce an overall Lipschitz constant of 1.

It is enlightening to view the proximal distance algorithm through the lens of concave-convex programming. Recall that the function
\begin{eqnarray}
s(\bx) & = & \sup_{\by \in C} \,\left[\by^t\bx - \frac{1}{2}\|\by\|^2\right]
\amp = \amp \frac{1}{2}\|\bx\|^2-\frac{1}{2}\dist(\bx, C)^2 \label{s_function}
\end{eqnarray}
is closed and convex for any nonempty closed set $C$. Danskin's theorem \citep{lange2016MM} justifies the directional derivative expression
\begin{eqnarray*}
d_{\bv}s(\bx) & = & \sup_{\by \in P_C(\bx)} \by^t\bv
\amp = \amp  \sup_{\by \in \conv P_C(\bx)} \by^t\bv.
\end{eqnarray*}
This equality allows us to identify the subdifferential 
$\partial s(\bx)$ as the convex hull $\conv P_C(\bx)$. For any $\by \in \partial s(\bx_k)$, the supporting hyperplane inequality entails 
\begin{eqnarray*}
\frac{1}{2}\dist(\bx, C)^2 & = & \frac{1}{2}\|\bx\|^2 -s(\bx) \\
& \le & \frac{1}{2}\|\bx\|^2-s(\bx_k)-\by^t(\bx-\bx_k) \\
& = & \frac{1}{2}\|\bx-\by\|^2 +d,
\end{eqnarray*}
where $d$ is a constant not depending on $\bx$. The same majorization can be generated by rearranging the majorization 
\begin{eqnarray*}
\frac{1}{2}\dist(\bx, C)^2 & \le & \frac{1}{2}\sum_i \beta_i \|\bx-\bp_i\|^2
\end{eqnarray*}
when $\by$ is the convex combination $\sum_i \beta_i \bp_i$ of vectors
$\bp_i$ from $P_C(\bx_k)$. These facts demonstrate that the proximal distance algorithm minimizing 
\begin{eqnarray*}
f(\bx)+\frac{\rho}{2}\dist(\bx, C)^2 
& = & f(\bx)+\frac{\rho}{2}\|\bx\|^2- \rho s(\bx)
\end{eqnarray*} 
is a special case of concave-convex programming when $f(\bx)$ is convex. It is worth emphasizing that $f(\bx)+\frac{\rho}{2}\|\bx\|^2$ is often strongly convex regardless of whether
$f(\bx)$ itself is convex. If we replace the penalty $\dist(\bx,C)^2$ by the penalty
$\dist(\bD\bx,C)^2$ for a matrix $\bD$, then the function $s(\bD\bx)$ is still closed and convex, and minimization of $f(\bx)+\frac{\rho}{2}\dist(\bD\bx, C)^2$ can also be viewed as an exercise in concave-convex programming.

In the presence of convexity, the proximal distance algorithm is guaranteed to converge. Our exposition relies on well-known operator results \citep{bauschke2011convex}. Proximal operators in general and projection operators in particular are nonexpansive and averaged. By definition an averaged operator 
\begin{eqnarray*}
M(\bx)  & = &  \alpha \bx+(1-\alpha) N(\bx)
\end{eqnarray*}
is a  convex combination of a nonexpansive 
operator $N(\bx)$ and the identity operator $\bI$. The averaged operators 
on $\mathbb{R}^p$ with $\alpha \in (0,1)$ form a convex set closed 
under functional composition. Furthermore, $M(\bx)$ and the base operator $N(\bx)$ share their fixed points. The celebrated theorem of \cite{krasnosel1955two} and \cite{mann1953mean} says that if an averaged operator $M(\bx) = \alpha \bx+(1-\alpha)N(\bx)$ possesses one or more fixed points, then the iteration scheme $\bx_{k+1} = M(\bx_k)$ converges to a fixed point. 

These results immediately apply to minimization of the penalized loss 
\begin{eqnarray}
h_\rho(\bx) & = &   f(\bx)+\frac{\rho}{2}\sum_{i=1}^m \alpha_i \dist(\bx,C_i)^2 .
\label{h_function}
\end{eqnarray}
Given the choice $\by_k = \sum_{i=1}^m \alpha_i P_{C_i}(\bx_k)$,
the algorithm map $\bx_{k+1} = \prox_{\rho^{-1}f}(\by_k)$ is an averaged operator, being the composition of two averaged operators. 
Hence, the Krasnosel'skii-Mann theorem guarantees convergence to a
fixed point if one or more exist.  Now $\bz$ is a fixed point if and only if
\begin{eqnarray*}
h_\rho(\bz) & \le &
f(\bx)+\frac{\rho}{2}\sum_{i=1}^m \alpha_i \|\bx-P_{C_i}(\bz)\|^2
\end{eqnarray*} 
for all $\bx$. In the presence of convexity, this is equivalent to
the directional derivative inequality
\begin{eqnarray*}
0 & \le & d_{\bv}f(\bz)+\rho \sum_{i=1}^m \alpha_i [\bz-P_{C_i}(\bz)]^t\bv 
\amp =  \amp d_{\bv}h_\rho(\bz)
\end{eqnarray*} 
for all $\bv$, which is in turn equivalent to $\bz$ minimizing
$h_\rho(\bx)$. Hence, if $h_\rho(\bx)$ attains its minimum value, then the proximal distance iterates converge to a minimum point.

Convergence of the overall proximal distance algorithm is tied to
the convergence of the classical penalty method \citep{beltrami1970algorithmic}.
In our setting, the loss is $f(\bx)$, and the penalty is
$q(\bx) = \frac{1}{2}\sum_{i=1}^m \alpha_i \dist(\bx,C_i)^2$. Assuming the objective $ f(\bx)+\rho q(\bx)$ is coercive, the theory mandates that the solution path $\bx_\rho$ is bounded and any cluster point of the path attains the minimum value of $f(\bx)$ subject to the constraints. Furthermore, if $f(\bx)$ is coercive and possesses a unique minimum point in the constraint set $C$, then the path $\bx_\rho$ converges to that point.

Proximal distance algorithms converge at a painfully slow rate.
Following \cite{mairal2013optimization}, one can readily exhibit a precise bound.  In the convex setting, we first observe that the surrogate function $g_\rho(\bx \mid \bx_k)$ is $\rho$-strongly convex. 
Consequently, the stationarity condition 
${\bf 0} \in \partial g_\rho(\bx_{k+1} \mid \bx_k)$ implies
\begin{eqnarray}
g_\rho(\bx \mid \bx_k) & \ge & g_\rho(\bx_{k+1} \mid \bx_k)
+ \frac{\rho}{2}\|\bx-\bx_{k+1}\|^2 \label{error_bd1}
\end{eqnarray}
for all $\bx$. In the notation (\ref{h_function}),  the difference
\begin{eqnarray*}
d_\rho(\bx \mid \bx_k) & = & g_\rho(\bx \mid \bx_k) -h_\rho(\bx)
\amp = \amp \frac{\rho}{2}\|\bx-\by_k\|^2 -\frac{\rho}{2} \sum_{i=1}^m \alpha_i \dist(\bx,C_i)^2 
\end{eqnarray*}
has a $\rho$-Lipschitz gradient because
\begin{eqnarray*}
\nabla d_\rho(\bx \mid \bx_k) & = & \rho(\bx-\by_k) -\rho \sum_{i=1}^m \alpha_i
[\bx-P_{C_i}(\bx)] \amp = \amp \rho \sum_{i=1}^m \alpha_i P_{C_i}(\bx)
- \rho \by_k.
\end{eqnarray*}
The tangency conditions $d_\rho(\bx_k \mid \bx_k)=0$ and 
$\nabla d_\rho(\bx_k \mid \bx_k)= {\bf 0}$ therefore yield 
\begin{eqnarray}
d_\rho(\bx \mid \bx_k ) & \le & d_\rho(\bx_k \mid \bx_k) +
 \nabla d_\rho(\bx_k)^t(\bx-\bx_k) +\frac{\rho}{2}\|\bx-\bx_{k}\|^2 \amp = \amp \frac{\rho}{2}\|\bx-\bx_{k}\|^2 \label{error_bd2}
\end{eqnarray}
for all $\bx$. At a minimum $\bz$ of $h_\rho(\bx)$, combining inequalities (\ref{error_bd1}) and (\ref{error_bd2}) gives
\begin{eqnarray*}
h_\rho(\bx_{k+1})+\frac{\rho}{2}\|\bz-\bx_{k+1}\|^2 & \le &
g_\rho(\bx_{k+1} \mid \bx_k)+\frac{\rho}{2}\|\bz-\bx_{k+1}\|^2 \\
& \le & g_\rho(\bz \mid \bx_k) \\
& = & h_\rho(\bz)+d_\rho(\bz \mid \bx_k) \\
& \le & h_\rho(\bz)+\frac{\rho}{2}\|\bz-\bx_{k}\|^2 .
\end{eqnarray*}
Adding the result 
\begin{eqnarray*}
h_\rho(\bx_{k+1})-h_\rho(\bz) & \le &
\frac{\rho}{2}\Big(\|\bz-\bx_{k}\|^2
 -\|\bz-\bx_{k+1}\|^2\Big)
\end{eqnarray*}
over $k$ and invoking the descent property $h_\rho(\bx_{k+1}) \le h_\rho(\bx_k)$ produce the desired error bound
\begin{eqnarray*}
h_\rho(\bx_{k+1})-h_\rho(\bz) & \le &
\frac{\rho}{2(k+1)} \Big(\|\bz-\bx_{0}\|^2 -\|\bz-\bx_{k+1}\|^2\Big)
\amp \le \amp \frac{\rho}{2(k+1)}\|\bz-\bx_{0}\|^2.
\end{eqnarray*}

The $O(\rho k^{-1})$ convergence rate of the proximal distance
algorithm suggests that one should slowly send $\rho$ to $\infty$ and refuse to wait until convergence occurs for any given $\rho$. It also suggests that Nesterov acceleration may vastly improve the chances for convergence. Nesterov acceleration for the general proximal gradient algorithm with loss $\ell(\bx)$ and penalty $p(\bx)$ takes the form
\begin{eqnarray}
\bz_k & = & \bx_k+\frac{k-1}{k+d-1}(\bx_k-\bx_{k-1}) \nonumber \\
\bx_{k+1} & = & \prox_{L^{-1}\ell}[\bz_k-L^{-1}\nabla p(\bz_k)], \label{Nesterov_acceleration}
\end{eqnarray}
where $L$ is the Lipschitz constant for $\nabla p(\bx)$ and $d$ is typically chosen to be 3. Nesterov acceleration achieves 
an $O(k^{-2})$ convergence rate \citep{su2014differential}, which is vastly superior to the $O(k^{-1})$ rate achieved by proximal gradient descent. The Nesterov update possesses the desirable property of preserving affine constraints. In other words, if $\bA \bx_{k-1} = \bb$ and $\bA\bx_k = \bb$, then $\bA \bz_k = \bb$ as well.  In 
subsequent examples, we will accelerate our proximal distance algorithms by applying the algorithm map $M(\bx)$ given by equation (\ref{algorithm_map}) to the shifted point $\bz_k$ of equation (\ref{Nesterov_acceleration}) , yielding the accelerated update $\bx_{k+1} = M(\bz_k)$.
Algorithm \ref{alg:proxdist} provides a schematic of a proximal distance algorithm with Nesterov acceleration.
The recent paper of Ghadimi and Lan \cite{Ghadimi2015} extends Nestorov acceleration to nonconvex settings. 

\begin{algorithm}
    \SetKwInOut{Input}{input}
    \SetKwInOut{Output}{output}

    \Input{}
        $\rho_{\mathrm{initial}} > 0$ \Comment*{an initial penalty value} 
        $\rho_{\mathrm{inc}} > 1$ \Comment*{the penalty increment} 
        $\rho_{\mathrm{max}}$ \Comment*{a maximum penalty value} 
        $K_{\mathrm{max}}$ \Comment*{the maximum number of iterations} 
        $k_{\rho}$ \Comment*{the increment frequency} 
        $f$ \Comment*{the function to optimize} 
        $P_C$ \Comment*{the projection onto the constraint set} 
        $\epsilon_{ \mathrm{loss} } > 0$ \Comment*{convergence tolerance for the loss function} 
        $\epsilon_{ \mathrm{dist} } > 0$ \Comment*{convergence tolerance for constraint feasibility} 
    \BlankLine
    \Output{ A vector $ \bx^{+} \approx \argmin_{\bx} f(\bx) $ 
    \Comment*{subject to the constraint $\bx \in C$}}
    \BlankLine
    $\rho \leftarrow \rho_{ \mathrm{initial} }$ \Comment*{Set initial penalty value}
    \BlankLine
    $q_{0} = q_{1} = \infty$ \Comment*{Track convergence of $f$}
    \BlankLine
    $d_{0} = d_{1} = \infty$ \Comment*{Track distance to constraint}
    \BlankLine
    $\bx_{0} = \bx_{1} = \bzero$ \Comment*{Set initial iterates to origin}
    \BlankLine
    \Comment{Main algorithm loop}
    \For{ $k = 2, \ldots, K_{\mathrm{max}}$ }{
        \BlankLine
        $\bz_{k} \leftarrow \bx_{ k - 1 } + \frac{ k - 1 }{ k + 2 } \left( \bx_{ k - 1 } - \bx_{ k - 2 } \right)$ \Comment*{Apply Nesterov acceleration}
        \BlankLine
        $\bx_{ k - 2 } \leftarrow \bx_{ k - 1 }$ \Comment*{Save penultimate iterate}
        \BlankLine
        $\bz_{ \mathrm{proj}, k } \leftarrow P_{C} \left( \bz_{k} \right)$ \Comment*{Project onto constraints}
        \BlankLine
        $\bx_{k} \leftarrow \prox_{ \rho f } \left( \bz_{ \mathrm{proj}, k } \right)$ \Comment*{Apply proximal distance update}
        \BlankLine
        $q_{k} \leftarrow f \left( \bx_{k} \right)$ \Comment*{Compute new loss}
        \BlankLine
        $d_{k} \leftarrow \left\| \bx_{k} - \bz_{ \mathrm{proj}, k } \right\|_{2}$ \Comment*{Compute new distance to $C$}
        \BlankLine
        \Comment{Exit if converged}
        \If{ $ \left| q_{k} - q_{ k - 1 } \right| < \epsilon_{ \mathrm{loss} }$ {\normalfont and} $ \left| d_{k} - d_{ k - 1 } \right| < \epsilon_{ \mathrm{dist} } $  }{
            \BlankLine
            \Return $\bx^{+} \leftarrow \bx_{k + 1}$ 
        }
        \Else{
            \BlankLine
            $q_{ k - 1 } \leftarrow q_{k}$ \Comment*{Save current loss}
            \BlankLine
            $d_{ k - 1 } \leftarrow d_{k}$ \Comment*{Save current distance to $C$}
            \BlankLine
            \Comment{Update penalty $\rho$ every $k_{\rho}$ iterations}
            \If{ $k = k_{\rho}$ }{

                $\rho \leftarrow \min\left( \rho_{ \mathrm{max} }, \rho \times \rho_{ \mathrm{inc} } \right)$

                $\bx_{ k - 1 } \leftarrow \bx_{k}$ \Comment*{Save previous iterate}
            }
        }
    }
    \caption{A typical proximal distance algorithm}
    \label{alg:proxdist}
\end{algorithm}

\section{Convergence: General Case}
\label{sec:nonconvex_convergence}

To simplify notation, we restrict attention to a single contraint
set $S$. Our strategy for addressing convergence relies on Zangwill's global convergence theorem \citep{luenberger1984linear}. This result depends in turn on the notion of a closed multi-valued map $\bN(\bx)$. If $\bx_k$ converges to $\bx_\infty$ and $\by_k \in N(\bx_k)$ converges to $\by_\infty$, then for $\bN(\bx)$ to be closed, we must have $\by_\infty \in N(\bx_\infty)$. The next proposition furnishes a prominent example.
\begin{proposition} \label{proposition1}
If $S$ is a closed nonempty set in $\mathbb{R}^p$, then the projection operator $P_S(\bx)$ is closed. Furthermore, if the sequence $\bx_k$ is bounded, then the set $\cup_k P_S(\bx_k)$ is bounded as well.
\end{proposition}
{\bf Proof:} Let $\bx_k$ converge to $\bx_\infty$ and $\by_k \in P_S(\bx_k)$ converge to $\by_\infty$. For an arbitrary $\by \in S$, taking limits in the inequality $\|\bx_k - \by_k\| \le \|\bx_k-\by\|$ yields $\|\bx_\infty - \by_\infty\| \le \|\bx_\infty-\by\|$; consequently, $\by_\infty \in P_S(\bx_\infty)$. To prove the second assertion, take $\by_k \in P_S(\bx_k)$ and observe that
\begin{eqnarray*}
\|\by_k\| & \le & \|\bx_k-\by_k\|+\|\bx_k\| \\
& \le & \|\bx_k-\by_1\|+ \|\bx_k\| \\
& \le & \|\bx_k-\bx_1\|+\|\bx_1-\by_1\|+ \|\bx_k\| \\
& \le & \|\bx_k\|+\|\bx_1\|+ \dist(\bx_1,S) + \|\bx_k\| ,
\end{eqnarray*}
which is bounded above by the constant $\dist(\bx_1,S) + 3\sup_{m \ge 1}\|\bx_m\|$.
\qed \svskip

Zangwill's global convergence theorem is phrased in terms of an algorithm map $M(\bx)$ and a real-valued objective $h(\bx)$. The theorem requires a critical set $\Gamma$ outside which $M(\bx)$ is closed. Furthermore, all iterates $\bx_{k+1} \in M(\bx_k)$ must fall within a compact set. Finally, the descent condition $h(\by) \le h(\bx)$ should hold for all $\by \in M(\bx)$, with strict inequality when $\bx \not\in \Gamma$. If these conditions are valid, then every convergent subsequence of $\bx_k$ tends to a point in $\Gamma$. In the proximal distance context, we define the complement of $\Gamma$ to consist of the points $\bx$ with 
\begin{eqnarray*}
f(\by)+\frac{\rho}{2}\dist(\by,S)^2 & < &
f(\bx)+\frac{\rho}{2}\dist(\bx,S)^2
\end{eqnarray*}
for all $\by \in M(\bx)$. This definition plus the monotonic nature of the proximal distance algorithm 
\begin{eqnarray*}
\bx_{k+1} & \in &  
M(\bx_k) \amp = \amp \bigcup_{\bz_k \in P_S(\bx_k)}\argmin_{\bx} \Big[f(\bx)+\frac{\rho}{2}\|\bx-\bz_k\|^2\Big] 
\end{eqnarray*}
force the satisfaction of Zangwill's final requirement. Note that if $f(\bx)$ is differentiable, then a point $\bx$ belongs to $\Gamma$ whenever 
${\bf 0 } \in \nabla f(\bx)+\rho \bx- \rho P_S(\bx)$. 

In general, the algorithm map $M(\bx)$ is multi-valued in two senses. First, for a given $\bz_k \in P_S(\bx_k)$, the minimum may be achieved at multiple points. This contingency is ruled out if the proximal map of $f(\bx)$ is unique. Second, because $S$ may be nonconvex, the projection may be multi-valued. This sounds distressing, but the points $\bx_k$ where this occurs are exceptionally rare. Accordingly, it makes no practical  difference that we restrict the anchor points $\bz_k$ to lie in $P_S(\bx_k)$ rather than 
in $\conv P_S(\bx_k)$.

\begin{proposition} \label{proposition2}
If $S$ is a closed nonempty set in $\mathbb{R}^p$, then the projection operator $P_S(\bx)$ is single valued except on a set of Lebesgue measure $0$.
\end{proposition}
{\bf Proof:} In fact, a much stronger result holds. Since the function $s(\bx)$ of equation (\ref{s_function}) is convex and finite,  Alexandrov's theorem \citep{niculescu2006convex} implies that it is almost everywhere twice differentiable.  In view of the identities $\frac{1}{2}\dist(\bx, S)^2 = \frac{1}{2}\|\bx\|^2 -s(\bx)$ and $\bx-P_S(\bx) = \nabla \frac{1}{2}\dist(\bx, S)^2$ where $P_S(\bx)$ is single valued, it follows that $P_S(\bx) = \nabla s(\bx)$ is almost everywhere differentiable. \qed

\begin{proposition} \label{proposition3}
The algorithm map $M(\bx)$ is everywhere closed.
\end{proposition}
{\bf Proof:}  If $\bx_k$ tends to $\bx_\infty$ and $\bz_k \in M(\bx_k)$ tends to $\bz_\infty$, then we must demonstrate that
$\bz_\infty \in M(\bx_\infty)$. By definition all $\bx$
satisfy
\begin{eqnarray}
f(\bz_{k})+\frac{\rho}{2}\|\bz_{k}-\by_k\|^2
& \le & f(\bx)+\frac{\rho}{2}\|\bx-\by_k\|^2
\label{map_descent}
\end{eqnarray}
for any $\by_k \in P_S(\bx_k)$. A sequence $\by_k$ of such 
values is bounded and therefore has a convergent subsequence with limit $\by_\infty$. 
Taking limits in inequality (\ref{map_descent}) along the subsequence
gives
\begin{eqnarray*}
f(\bz_\infty)+\frac{\rho}{2}\|\bz_\infty-\by_\infty \|^2
& \le & f(\bx)+\frac{\rho}{2}\|\bx-\by_\infty\|^2.
\end{eqnarray*} 
Because $P_S(\bx)$ is a closed map, $\by_\infty \in P_S(\bx_\infty)$ and consequently $\bz_\infty \in M(\bx_\infty)$. \qed \svskip

To apply Zangwill's global convergence theory, we must in addition
prove that the iterates $\bx_{k+1} = M(\bx_k)$ remain within a compact
set. This is true whenever the objective is coercive since the 
algorithm is a descent algorithm. As noted earlier, the coercivity
of $f(\bx)$ is a sufficient condition. One can readily concoct
other sufficient conditions. For example, if $f(\bx)$ is bounded below,
say nonnegative, and $S$ is compact, then the objective is
also coercive. Indeed, if $S$ is contained in the ball of radius $r$
about the origin, then
\begin{eqnarray*}
\|\bx\| & \le & \|\bx-P_S(\bx)\|+\|P_S(\bx)\|
\amp \le \amp \dist(\bx,S)+r ,
\end{eqnarray*}
which proves that $\dist(\bx,S)$ is coercive.

\begin{proposition} \label{proposition4}
If $S$ is closed and nonempty, the objective $f(\bx)+\frac{1}{2}\dist(\bx,S)^2$ is coercive, and the proximal operator $\prox_{\rho^{-1} f}(\bx)$ is everywhere nonempty,
then all limit points of the iterates $\bx_{k+1} \in M(\bx_k)$ of the proximal distance algorithm occur in the critical set $\Gamma$.
\end{proposition}
{\bf Proof:} See the foregoing discussion. \qed \svskip

This result is slightly disappointing. A limit point $\bx$ could potentially exist with improvement in the objective for some but not all $\by \in \conv P_S(\bx)$. This fault
is mitigated by the fact that $P_S(\bx)$ is almost always single valued. In common with 
other algorithms in nonconvex optimization, we also cannot rule out convergence to a local minimum or a saddlepoint. 

One can improve on Proposition \ref{proposition4} by assuming that the surrogates 
$g_\rho(\bx \mid \bx_k)$ are all $\mu$-strongly convex. This is a small concession to make because $\rho$ is typically large.
\begin{proposition} \label{proposition5}
Under the $\mu$-strongly convexity assumption on the surrogates $g_\rho(\bx \mid \bx_k)$,
the proximal distance iterates satisfy $\lim_{k \to \infty} \|\bx_{k+1}-\bx_k\| = 0$.
As a consequence, the set of limit points is closed and connected. Furthermore,
if each limit point is isolated, then the iterates converge to a critical point.
\end{proposition}
{\bf Proof:}  The strong-convexity inequality 
\begin{eqnarray*}
g_\rho(\bx_k \mid \bx_k) & \ge &
g_\rho(\bx_{k+1} \mid \bx_k)+\frac{\mu}{2}\|\bx_k-\bx_{k+1}\|^2 
\end{eqnarray*}
and the tangency and domination properties of the algorithm imply
\begin{eqnarray}
h_\rho(\bx_k)-h_\rho(\bx_{k+1}) & \ge & \frac{\mu}{2}\|\bx_k-\bx_{k+1}\|^2. 
\label{strong_convexity_ineq}
\end{eqnarray}
Since the difference in function values tends to $0$, this validates the stated limit. The remaining assertions follow from Propositions 7.3.3 and 7.3.5 of \citep{lange2016MM}.
\qed \svskip

Further progress requires even more structure. Fortunately, what we now pursue 
applies to generic MM algorithms. We start with the concept of a Fr\'echet subdifferential \citep{kruger2003frechet}. If $h(\bx)$ is a function mapping $\mathbb{R}^p$ into $\mathbb{R} \cup \{+\infty\}$, then its Fr\'echet subdifferential at 
$\bx \in \dom f$ is the set
\begin{eqnarray*}
\partial^F h(\bx) & = & \Big\{\bv: \liminf_{\by \to \bx}
\frac{h(\by)-h(\bx) - \bv^t(\by-\bx)}{\|\by-\bx\|} \ge 0 \Big\}.
\end{eqnarray*}
The set $\partial^F h(\bx)$ is closed, convex, and possibly empty. If $h(\bx)$ is convex, then $\partial^F h(\bx)$ reduces to its convex subdifferential. If $h(\bx)$ is differentiable, then $\partial^F h(\bx)$ reduces to its ordinary differential. At a local minimum $\bx$, Fermat's rule ${\bf 0} \in \partial^F h(\bx)$ holds. 

\begin{proposition}
In an MM algorithm, suppose that $h(\bx)$ is coercive, that the surrogates $g(\bx \mid \bx_k)$ are differentiable, and that the algorithm map $M(\bx)$ is closed. Then every limit point $\bz$ of the MM sequence $\bx_k$ is critical in the sense that 
${\bf 0} \in \partial^F (-h)(\bz)$. 
\end{proposition}
{\bf Proof:} Let the subsequence $x_{k_m}$ of the MM sequence 
$\bx_{k+1} \in M(\bx_k)$ converge to $\bz$. By passing to a subsubsequence if necessary, we may suppose that $\bx_{k_m+1}$ converges to $\by$. Owing to our closedness assumption, $\by \in M(\bz)$. Given that $h(\by)=h(\bz)$, it is obvious that $\bz$ also minimizes $g(\bx \mid \bz)$ and that 
${\bf 0} = \nabla g(\bz \mid \bz)$.  Since the difference $\Delta(\bx \mid \bz)= g(\bx \mid \bz)-h(\bx)$ achieves its minimum at $\bx = \bz$, the Fr\'echet subdifferential $\partial^F \Delta(\bx \mid \bz)$ satisfies 
\begin{eqnarray*}
{\bf 0} & \in & \partial^F \Delta(\bz \mid \bz) \amp = \amp
\nabla g(\bz \mid \bz) + \partial^F (-h) (\bz).
\end{eqnarray*}
It follows that ${\bf 0} \in \partial^F (-h) (\bz)$.
\qed \svskip

We will also need to invoke Lojasiewicz's inequality. This deep result depends on some rather arcane algebraic geometry \citep{bierstone1988semianalytic,bochnak2013real}. It applies to semialgebraic functions and their more inclusive cousins semianalytic functions and subanalytic functions. For simplicity we focus on semialgebraic functions. The class of semialgebraic subsets of $\mathbb{R}^p$ is the smallest class that:
\begin{description}
\item[a)] contains all sets of the form
$\{\bx: q(\bx)>0 \}$ for a polynomial $q(\bx)$ in $p$ variables,
\item[b)] is closed under the formation of finite unions, finite intersections, and set complementation. 
\end{description}
A function $a:\mathbb{R}^p \mapsto \mathbb{R}^r$ is said to be semialgebraic if its graph is a semialgebraic set of $\mathbb{R}^{p+r}$. The class of real-valued semialgebraic functions contains all polynomials $p(\bx)$ and all $0$/$1$ indicators of algebraic sets.
It is closed under the formation of sums, products, absolute values,
reciprocals when $a(\bx) \ne 0$, $n$th roots when $a(\bx) \ge 0$, and maxima $\max\{a(\bx),b(\bx)\}$ and minima $\min\{a(\bx),b(\bx)\}$. For our purposes, it is important to note that $\dist(\bx, S)$ is a semialgebraic function whenever $S$ is a semialgebraic set.

Lojasiewicz's inequality in its modern form \citep{bolte2007lojasiewicz} requires a
function $h(\bx)$ to be closed (lower semicontinuous) and subanalytic with a closed domain. If $\bz$ is a critical point of $h(\bx)$, then 
\begin{eqnarray*}
|h(\bx) - h(\bz)|^{\theta} & \le & c \|\bv \|
\end{eqnarray*}
for all $\bx \in B_{r}(\bz) \cap \dom\partial^F h$ satisfying $h(\bx)>h(\bz)$ and all $\bv$ in $\partial^F h(\bx)$. Here the exponent $\theta \in [0,1)$, the radius $r$, and the constant $c$ 
depend on $\bz$. This inequality is valid for  semialgebraic functions since they are automatically subanalytic. We will apply Lojasiewicz's inequality to the limit points of an MM algorithm. The next proposition is an elaboration and expansion of known results \citep{attouch2010proximal,bolte2007lojasiewicz,cui2018composite,kang2015global,le2009convergence}.

\begin{proposition} In an MM algorithm suppose the objective $h(\bx)$ is coercive, continuous, and subanalytic and all surrogates $g(\bx \mid \bx_k)$ are continuous, 
$\mu$-strongly convex, and satisfy the Lipschitz condition
\begin{eqnarray*}
\|\nabla g(\ba \mid \bx_k)-\nabla g(\bb \mid \bx_k) \|
& \le & L \|\ba-\bb\| 
\end{eqnarray*}
on the compact set $\{\bx: h(\bx) \le h(\bx_0)\}$. Then the MM iterates 
$\bx_{k+1} = \argmin_{\bx} g(\bx \mid \bx_k)$ converge to a critical point.  
\end{proposition} 
{\bf Proof:} Because $\Delta(\bx \mid \by)= g(\bx \mid \by)-h(\bx)$ achieves its minimum at $\bx = \by$, the Fr\'echet subdifferential $\partial^F \Delta(\bx \mid \by)$ satisfies 
\begin{eqnarray*}
{\bf 0} & \in & \partial^F \Delta(\by \mid \by) \amp = \amp
\nabla g(\by \mid \by) + \partial^F (-h) (\by).
\end{eqnarray*}
It follows that $-\nabla g(\by \mid \by) \in \partial^F (-h) (\by)$.
Furthermore, by assumption
\begin{eqnarray*}
\|\nabla g(\ba \mid \bx_k)-\nabla g(\bb \mid \bx_k) \|
& \le & L\|\ba-\bb\| 
\end{eqnarray*}
for all relevant $\ba$ and $\bb$ and $\bx_k$. In particular, because $\nabla g(\bx_{k+1} \mid \bx_k) = {\bf 0}$, we have
\begin{eqnarray}
\|\nabla g(\bx_k \mid \bx_k) \|
& \le & L \|\bx_{k+1}-\bx_k\|. \label{lipschitz_MM}
\end{eqnarray}

Let $W$ denote the set of limit points. The objective $h(\bx)$ is constant on $W$ with value $\bar{h} = \lim_{k \to \infty} h(\bx_k)$. According to the Lojasiewicz inequality applied for the subanalytic function $\bar{h}-h(\bx)$, for each $\bz \in W$ there exists an open ball $B_{r(\bz)}(\bz)$ of radius $r(\bz)$ around $\bz$ and an exponent 
$\theta(\bz) \in [0,1)$ such that
\begin{eqnarray*}
|h(\bu) - h(\bz)|^{\theta(\bz)} & = & 
|\bar{h} - h(\bu) - \bar{h}+\bar{h}|^{\theta(\bz)}
\amp \le \amp c(\bz) \|\bv \|
\end{eqnarray*}
for all $\bu \in B_{r(\bz)}(\bz)$ and all 
$\bv \in \partial^F (\bar{h}-h) (\bu)= \partial^F (-h) (\bu)$. We will apply this inequality to $\bu=\bx_k$ and 
$\bv = -\nabla g(\bx_k \mid \bx_k)$.
In so doing, we would like to assume that the exponent $\theta(\bz)$ and constant $c(\bz)$ do not depend on $\bz$. With this end in mind, cover the compact set  
$W$ by a finite number of balls $B_{r(\bz_i)}(\bz_i)$ and take 
$\theta=\max_i \theta(\bz_i)<1$ and $c = \max_i c(\bz_i)$. For a sufficiently large $K$, every $\bx_k$ with $k \ge K$ falls within one
of these balls and satisfies $|\bar{h}-h(\bx_k)| < 1$. Without loss
of generality assume $K=0$. The Lojasiewicz inequality reads
\begin{eqnarray}
|\bar{h}-h(\bx_k)|^{\theta} & \le & c \|\nabla g(\bx_k \mid \bx_k) \|. 
\label{Lojasiewicz_MM}
\end{eqnarray}
In combination with the concavity of the function $t^{1-\theta}$ on $[0,\infty)$, inequalities (\ref{strong_convexity_ineq}),
(\ref{lipschitz_MM}), and (\ref{Lojasiewicz_MM}) imply
\begin{eqnarray*}
[h(\bx_k) - \bar{h}]^{1-\theta}-[h(\bx_{k+1})-\bar{h}]^{1-\theta} & \ge & 
\frac{1-\theta}{[h(\bx_k)-\bar{h}]^{\theta}}[h(\bx_k)-h(\bx_{k+1})] \\
& \ge & \frac{1-\theta}{c\|\nabla g(\bx_k \mid \bx_k) \|}
\frac{\mu}{2}\|\bx_{k+1}-\bx_k\|^2 \\
& \ge & \frac{(1-\theta)\mu}{2cL}\|\bx_{k+1}-\bx_k\|.
\end{eqnarray*}
Rearranging this inequality and summing over $k$ yield
\begin{eqnarray*}
\sum_{n=0}^\infty \|\bx_{k+1}-\bx_k\| & \le &
\frac{2cL}{(1-\theta)\mu}[h(\bx_0) - \bar{h}]^{1-\theta}
\end{eqnarray*}
Thus, the sequence $\bx_k$ is a fast Cauchy sequence and converges to a unique limit in $W$.
\qed \svskip

The last proposition applies to proximal distance algorithms. The loss $f(\bx)$ 
must be subanalytic and differentiable with a locally Lipschitz gradient. Furthermore,
all surrogates $g(\bx \mid \bx_k)=f(\bx)+\frac{\rho}{2}\|\bx-\by_k\|^2$ should be
coercive and $\mu$-strongly convex. Finally, the constraints sets $S_i$ should be
subanalytic. Semialgebraic sets and functions will do. Under these conditions and
regardless of how the projected points $P_{S_i}(\bx)$ are chosen, the MM iterates are guaranteed to converge to a critical point.

\section{Examples}
\label{sec:examples}

The following examples highlight the versatility of proximal distance algorithms
in a variety of convex and nonconvex settings. Programming details matter
in solving these problems. Individual programs are not necessarily long, but 
care must be exercised in projecting onto constraints, choosing tuning
schedules, folding constraints into the domain of the loss, implementing
acceleration, and declaring convergence. 
All of our examples are coded in the Julia programming language.
Whenever possible, competing software was run in the Julia environment via
the Julia module \texttt{MathProgBase} \citep{DunningHuchetteLubin2015,LubinDunningIJOC}. 
The sparse PCA problem relies on the software of Witten et al. \citep{wittentibshiranihastie2009}, 
which is coded in R. 
Convergence is tested at iteration $k$ by the two criteria
\begin{eqnarray*}
|f(\bx_k)- f(\bx_{k-1})| & \le & \epsilon_1 [|f(\bx_{k-1})|+1]  \quad \text{and} \quad 
\dist(\bx_k,C) \amp \le \amp \epsilon_2,
\end{eqnarray*}
where $\epsilon_1 = 10^{-6}$ and $\epsilon_2 = 10^{-4}$ are typical
values. The number of iterations until convergence is about 1000 in
most examples. This handicap is offset by the simplicity of each
stereotyped update. Our code is available as supplementary
material to this paper. Readers are encouraged to try the code and adapt
it to their own examples. 

\subsection{Linear Programming}
\label{ex:lp}

Two different tactics suggest themselves for constructing a proximal distance algorithm. The first tactic rolls the standard affine constraints $\bA \bx = \bb$ into the domain of the loss function $\bv^t\bx$. The standard nonnegativity requirement $\bx \ge {\bf 0}$ is achieved by penalization.
Let $\bx_k$ be the current iterate and $\by_k = (\bx_k)_{+}$ be its projection onto $\mathbb{R}^n_+$. Derivation of the proximal distance algorithm relies on the Lagrangian
\begin{eqnarray*}
\bv^t\bx+\frac{\rho}{2}\|\bx- \by_k\|^2+\blambda^t(\bA\bx-\bb).
\end{eqnarray*}
One can multiply the corresponding stationarity equation
\begin{eqnarray*}
{\bf 0} & = & \bv+\rho(\bx-\by_k)+\bA^t\blambda
\end{eqnarray*}
by $\bA$ and solve for the Lagrange multiplier $\blambda$ in the form
\begin{eqnarray}
\blambda & = & (\bA\bA^t)^{-1}(\rho \bA\by_k-\rho \bb -\bA \bv), 
\label{eq:lp_lambda_update}
\end{eqnarray}
assuming $\bA$ has full row rank. Inserting this value into the stationarity equation gives the MM update
\begin{eqnarray}
\label{eq:lp_update}
\bx_{k+1} & = & \by_k - \frac{1}{\rho}\bv
-\bA^{-}\left( \bA\by_k-\bb-\frac{1}{\rho} \bA \bv\right),
\end{eqnarray}
where $\bA^{-}= \bA^t (\bA\bA^t)^{-1}$ is the pseudo-inverse of $\bA$.

The second tactic folds the nonnegativity constraints into the
domain of the loss. Let $\bp_k$ denote the projection of $\bx_k$
onto the affine constraint set $\bA\bx=\bb$. Fortunately, the surrogate function $\bv^t\bx+\frac{\rho}{2}\|\bx-\bp_k\|^2$ splits the parameters.
Minimizing one component at a time gives the update $\bx_{k+1}$
with components
\begin{eqnarray}
x_{k+1,j} & = & \max\Big\{p_{kj} -\frac{v_j}{\rho},0\Big\}.
\label{eq:2nd_lp_update}
\end{eqnarray}
The projection $\bp_k$ can be computed via 
\begin{eqnarray}
\bp_k & = & \bx_k - \bA^{-}(\bA\bx_k-\bb),
\label{2nd_lp_projection}
\end{eqnarray}
where $\bA^{-}$ is again the pseudo-inverse of $\bA$. 

Table \ref{tab:lp} compares the accelerated versions of these two proximal distance algorithms to two efficient solvers. 
The first is the open-source Splitting Cone Solver (SCS) \citep{o2013conic}, which relies on a fast implementation of ADMM.
The second is the commercial Gurobi solver, which ships with implementations of both the simplex method and a barrier (interior point) method; in this example, we use its barrier algorithm. 
The first seven rows of the table summarize linear programs with dense data $\bA$, $\bb$, and $\bv$. The bottom six rows rely on random sparse matrices $\bA$ with sparsity level $0.01$. For dense problems, the proximal distance algorithms start the penalty constant  $\rho$ at 1 and double it every 100 iterations. Because we precompute and cache the pseudoinverse $\bA^{-}$ of $\bA$, the updates \eqref{eq:lp_update} and (\ref{eq:2nd_lp_update}) reduce to vector additions and matrix-vector multiplications. 

For sparse problems the proximal distance algorithms update $\rho$ by a factor of 1.5 every 50 iterations. To avoid computing large pseudoinverses, we appeal to the LSQR variant of the conjugate gradient method \citep{paigesaunders1982b,paigesaunders1982a} to solve the linear systems \eqref{eq:lp_lambda_update} and (\ref{2nd_lp_projection}). The optima of all four methods agree to about 4 digits of accuracy.  It is hard to declare an absolute winner in these comparisons. Gurobi and SCS clearly perform better on low-dimensional problems, but the proximal distance algorithms are competitive as dimensions increase. PD1, the proximal distance algorithm over an affine domain, tends to be more accurate than PD2. If high accuracy is not a concern, then the proximal distance algorithms are easily accelerated with a more aggressive update schedule for $\rho$. 

\begin{table}
	\centering
	\begin{tabular}{cccccccccc}
		\toprule
		\multicolumn{2}{c}{Dimensions} & \multicolumn{4}{c}{Optima} & \multicolumn{4}{c}{CPU Times (secs)} \\
		\cmidrule(r){1-2} \cmidrule(r){3-6} \cmidrule(r){7-10}
		$m$ & $n$ & PD1 & PD2 & SCS & Gurobi & PD1 & PD2 & SCS & Gurobi \\
		\hline
		2 & 4 & 0.2629 & 0.2629 & 0.2629 & 0.2629 & 0.0142 & 0.0010 & 0.0034 & 0.0038 \\
		4 & 8 & 1.0455 & 1.0457 & 1.0456 & 1.0455 & 0.0212 & 0.0021 & 0.0009 & 0.0011 \\
		8 & 16 & 2.4513 & 2.4515 & 2.4514 & 2.4513 & 0.0361 & 0.0048 & 0.0018 & 0.0029 \\
		16 & 32 & 3.4226 & 3.4231 & 3.4225 & 3.4223 & 0.0847 & 0.0104 & 0.0090 & 0.0036 \\
		32 & 64 & 6.2398 & 6.2407 & 6.2397 & 6.2398 & 0.1428 & 0.0151 & 0.0140 & 0.0055 \\
		64 & 128 & 14.671 & 14.674 & 14.671 & 14.671 & 0.2117 & 0.0282 & 0.0587 & 0.0088 \\
		128 & 256 & 27.116 & 27.125 & 27.116 & 27.116 & 0.3993 & 0.0728 & 0.8436 & 0.0335 \\
		256 & 512 & 58.501 & 58.512 & 58.494 & 58.494 & 0.7426 & 0.1538 & 2.5409 & 0.1954 \\
		512 & 1024 & 135.35 & 135.37 & 135.34 & 135.34 & 1.6413 & 0.5799 & 5.0648 & 1.7179 \\
		1024 & 2048 & 254.50 & 254.55 & 254.47 & 254.48 & 2.9541 & 3.2127 & 3.9433 & 0.6787 \\
		2048 & 4096 & 533.29 & 533.35 & 533.23 & 533.23 & 7.3669 & 17.318 & 25.614 & 5.2475 \\
		4096 & 8192 & 991.78 & 991.88 & 991.67 & 991.67 & 30.799 & 95.974 & 98.347 & 46.957 \\
		8192 & 16384 & 2058.8 & 2059.1 & 2058.5 & 2058.5 & 316.44 & 623.42 & 454.23 & 400.59 \\
		\bottomrule
	\end{tabular}
	\caption{CPU times and optima for linear programming. Here $m$ is the number of constraints, $n$ is the number of variables, PD1 is the proximal distance algorithm over an affine domain, PD2 is the proximal distance algorithm over a nonnegative domain, SCS is the Splitting Cone Solver, and Gurobi is the Gurobi solver. After $m = 512$ the constraint matrix $\boldsymbol{A}$ is initialized to be sparse with sparsity level $s = 0.01$.}
	\label{tab:lp}
\end{table}

\subsection{Constrained Least Squares \label{ex:nqp}}

Constrained least squares programming subsumes constrained quadratic programming.
A typical quadratic program involves minimizing the quadratic $\frac{1}{2} \bx^t \bQ \bx - \bp^t \bx$ subject to 
$\bx \in C$ for a positive definite matrix $\bQ$. Quadratic programming can be reformulated as least squares by taking the Cholesky decomposition $\bQ=\bL\bL^t$ of 
$\bQ$ and noting that
\begin{eqnarray*}
\frac{1}{2} \bx^t \bQ \bx - \bp^t \bx & = & 
\frac{1}{2}\|\bL^{-1}\bp-\bL^t\bx\|^2 - \frac{1}{2}\|\bL^{-1}\bp\|^2.
\end{eqnarray*}
The constraint $\bx \in C$ applies in both settings. It is particularly advantageous to reframe a quadratic program as a least squares problem when $\bQ$ is already presented in factored form or when it is nearly 
singular \citep{bemporad2018semidefiniteQP}. To simplify subsequent notation, we replace $\bL^t$ by the rectangular matrix $\bA$ and $\bL^{-1}\bp$ by $\by$. The key to solving constrained least squares is to express the proximal distance surrogate as
\begin{eqnarray*}
\frac{1}{2}\|\by-\bA\bx\|^2+\frac{\rho}{2}\|\bx - P_C(\bx_k)\|^2
& = & \frac{1}{2}
\left\| \begin{pmatrix} \by
\\ \sqrt{\rho}P_C(\bx_k) \end{pmatrix} - \begin{pmatrix} \bA \\
\sqrt{\rho}\bI \end{pmatrix}\bx
\right\|^2
\end{eqnarray*}
as in equation (\ref{johann_credit}). As noted earlier, in sparse problems the
update $\bx_{k+1}$ can be found by a fast stable conjugate gradient solver.

\begin{table}
	\centering
	\begin{tabular}{cccccccc}
			\toprule
			\multicolumn{2}{c}{Dimensions} & \multicolumn{3}{c}{Optima} & \multicolumn{3}{c}{CPU Times} \\
			\cmidrule(r){1-2} \cmidrule(r){3-5} \cmidrule(r){6-8}
            $n$ & $p$ & PD & IPOPT & Gurobi & PD & IPOPT & Gurobi \\
			\hline
             16    &    8 & 4.1515 & 4.1515 & 4.1515 & 0.0038 & 0.0044 & 0.0010\\
             32    &   16 & 10.8225 & 10.8225 & 10.8225 & 0.0036 & 0.0039 & 0.0010\\
             64    &   32 & 29.6218 & 29.6218 & 29.6218 & 0.0079 & 0.0079 & 0.0019\\
             128   &   64 & 43.2626 & 43.2626 & 43.2626 & 0.0101 & 0.0078 & 0.0033\\
             256   &  128 & 111.7642 & 111.7642 & 111.7642 & 0.0872 & 0.0151 & 0.0136\\
             512   &  256 & 231.6455 & 231.6454 & 231.6454 & 0.1119 & 0.0710 & 0.0619\\
             1024  &  512 & 502.1276 & 502.1276 & 502.1276 & 0.2278 & 0.4013 & 0.2415\\
             2048  & 1024 & 994.2447 & 994.2447 & 994.2447 & 1.2575 & 2.3346 & 1.1682\\
             4096  & 2048 & 2056.8381 & 2056.8381 & 2056.8381 & 1.3253 & 15.2214 & 7.4971\\
             8192  & 4096 & 4103.4611 & 4103.4611 & 4103.4611 & 3.0289 & 146.1604 & 49.7411\\
             16384 & 8192 & 8295.2136 & 8295.2136 & 8295.2136 & 6.8739 & 732.1039 & 412.3612\\
			\bottomrule
		\end{tabular}
	\caption{CPU times and optima for simplex-constrained least squares. Here $\boldsymbol{A} \in \mathbb{R}^{n \times p}$, PD is the proximal distance algorithm, IPOPT is the Ipopt solver, and Gurobi is the Gurobi solver. After $n = 1024$, the predictor matrix $\boldsymbol{A}$ is sparse.}
	\label{tab:nqp}
\end{table}

Table \ref{tab:nqp} compares the performance of the proximal distance algorithm for least squares estimation with probability-simplex constraints to the open source nonlinear interior point solver Ipopt \citep{waechterbiegler2005,waechterbiegler2006} and the interior point method of Gurobi. Simplex constrained problems arise in
hyperspectral imaging \citep{HeylenBurazerovicScheunders11SimplexLS,keshava2003survey}, portfolio optimization \citep{markowitz1952portfolio}, and density estimation \citep{bunea2010spades}.
Test problems were generated by filling an $n \times p$ matrix $\bA$ and an $n$-vector $\by$ with standard normal deviates. For sparse problems we set the sparsity level of $\bA$ to be $10/p$. Our setup ensures that $\bA$ has full rank
and that the quadratic program has a solution.
For the proximal distance algorithm, we start $\rho$ at 1 and multiply it by $1.5$ every 200 iterations. Table \ref{tab:nqp} suggests that the proximal distance algorithm and the interior point solvers perform equally well on small dense problems. However, in high-dimensional and low-accuracy environments, the proximal distance algorithm provides much better scalability. 

\subsection{Closest Kinship Matrix}
\label{ex:kin}

In genetics studies, kinship is measured by the fraction of genes two individuals share
identical by descent. For a given pedigree, the kinship coefficients for all pairs of individuals appear as entries in a symmetric kinship matrix $\bY$.  This matrix possesses 
three crucial properties: a) it is positive semidefinite, b) its entries are nonnegative, 
and c) its diagonal entries are $\frac{1}{2}$ unless some pedigree members are inbred.
Inbreeding is the exception rather than the rule.  Kinship matrices can be 
estimated empirically from SNP (single nucleotide polymorphisms) data,
but there is no guarantee that the three highlighted properties are satisfied. Hence,
it helpful to project $\bY$ to the nearest qualifying matrix. 

This projection problem is best solved by folding the positive semidefinite
constraint into the domain of the Frobenius loss function $\frac{1}{2}\|\bX-\bY\|_F^2$.
As we shall see, the alternative of imposing two penalties rather than one 
is slower and less accurate. Projection onto the constraints implied by conditions b) and c) is trivial. All diagonal entries $x_{ii}$ of $\bX$ are reset to $\frac{1}{2}$, and all off-diagonal entries $x_{ij}$ are reset to $\max\{x_{ij},0\}$. If $P(\bX_k)$ denotes the current projection, then the proximal distance algorithm minimizes the surrogate
\begin{eqnarray*}
g(\bX \mid \bX_{k}) & = & \frac{1}{2}\|\bX-\bY\|_F^2 + \frac{\rho}{2}\|\bX-P(\bX_k)\|_F^2 \\
& = & \frac{1+\rho}{2}\left\|\bX-\frac{1}{1+\rho}\bY - \frac{\rho}{1+\rho}P(\bX_k)\right\|_F^2 +c_k,
\end{eqnarray*}
where $c_k$ is an irrelevant constant. The minimum is found by extracting the
spectral decomposition $\bU \bD \bU^t$ of 
$\frac{1}{1+\rho}\bY + \frac{\rho}{1+\rho}P(\bX_k)$ and truncating the negative
eigenvalues. This gives the update $\bX_{k+1} = \bU \bD_+ \bU^t$ in obvious
notation.  This proximal distance algorithm and its Nesterov acceleration are simple 
to implement in a numerically oriented language such as Julia. The most onerous 
part of the calculation is clearly the repeated eigen-decompositions. 

Table \ref{tab:kin} compares three versions of the proximal distance algorithm
to Dykstra's algorithm \citep{boyle1986method}. Higham proposed Dykstra's algorithm for the related problem of finding the closest correlation matrix \cite{higham2002computing}.
In Table \ref{tab:kin} algorithm PD1 is the unadorned proximal distance algorithm, PD2 is the accelerated proximal distance, and PD3 is the accelerated proximal distance algorithm with the positive semidefinite constraints folded into the domain of the loss. On this
demanding problem, these algorithms are comparable to Dykstra's algorithm in speed but slightly less accurate. 
Acceleration of the proximal distance algorithm is effective in reducing both execution time and error. 
Folding the positive semidefinite constraint into the domain of the loss function leads to further improvements. 
The data matrices $\bM$ in these trials were populated by standard normal deviates and then symmetrized by averaging opposing off-diagonal entries. In algorithm PD1 we set $\rho_k = \max\{1.2^k,2^{22}\}$. In the accelerated versions PD2 and PD3 we started $\rho$ at 1 and multiplied it by 5 every 100 iterations. At the expense of longer compute times, better accuracy can be achieved by all three proximal distance algorithms with a less aggressive update schedule. 
\begin{table}
    \centering
    \begin{tabular}{crrrrrrrr}
        \toprule
        \multicolumn{1}{c}{Size} & \multicolumn{2}{c}{PD1} & \multicolumn{2}{c}{PD2}   & \multicolumn{2}{c}{PD3} & \multicolumn{2}{c}{Dykstra} \\
        \cmidrule(r){1-1} \cmidrule(r){2-3} \cmidrule(r){4-5} \cmidrule(r){6-7}
        \cmidrule(r){8-9}
        \multicolumn{1}{c}{$n$}  & \multicolumn{1}{r}{Loss} & \multicolumn{1}{r}{Time} & 
        \multicolumn{1}{r}{Loss} & \multicolumn{1}{r}{Time} & \multicolumn{1}{r}{Loss} & 
        \multicolumn{1}{r}{Time} & \multicolumn{1}{r}{Loss} & \multicolumn{1}{r}{Time} \\
        \hline
          2 & 1.64 & 0.36 & 1.64 & 0.01 & 1.64 & 0.01 & 1.64 & 0.00 \\
          4 & 2.86 & 0.10 & 2.86 & 0.01 & 2.86 & 0.01 & 2.86 & 0.00 \\
          8 & 18.77 & 0.21 & 18.78 & 0.03 & 18.78 & 0.03 & 18.78 & 0.00 \\
         16 & 45.10 & 0.84 & 45.12 & 0.18 & 45.12 & 0.12 & 45.12 & 0.02 \\
         32 & 169.58 & 4.36 & 169.70 & 0.61 & 169.70 & 0.52 & 169.70 & 0.37 \\
         64 & 837.85 & 16.77 & 838.44 & 2.90 & 838.43 & 2.63 & 838.42 & 4.32 \\
        128 & 3276.41 & 91.94 & 3279.44 & 18.00 & 3279.25 & 14.83 & 3279.23 & 19.73 \\
        256 & 14029.07 & 403.59 & 14045.30 & 89.58 & 14043.59 & 64.89 & 14043.46 & 72.79 \\
        \bottomrule
    \end{tabular}
    \caption{CPU times and optima for the closest kinship matrix problem. Here the kinship matrix is $n \times n$, PD1 is the proximal distance algorithm, PD2 is the accelerated proximal distance, PD3 is the accelerated proximal distance algorithm with the 
    positive semidefinite constraints folded into the domain of the loss, and Dykstra
    is Dykstra's adaptation of alternating projections. All times are in seconds.}
    \label{tab:kin}
\end{table}

\subsection{Projection onto a Second-Order Cone Constraint}
\label{ex:socp}

Second-order cone programming is one of the unifying themes of
convex analysis \citep{alizadeh2003second,lobo1998applications}. It
revolves around conic constraints of the form $\{\bu: \|\bA\bu+\bb\| \le \bc^t\bu+d\}$.
Projection of a vector $\bx$ onto such a constraint is facilitated 
by parameter splitting. In this setting parameter splitting introduces a vector 
$\bw$, a scalar $r$, and the two affine constraints $\bw=\bA\bu+\bb$ and
$r = \bc^t\bu+d$.  The conic constraint then reduces to the Lorentz cone
constraint $\|\bw\| \le r$, for which projection is straightforward
\citep{boyd2009convex}.  If we concatenate the parameters into the single
vector
\begin{eqnarray*}
\by & = & \begin{pmatrix} \bu \\ \bw \\ r \end{pmatrix}
\end{eqnarray*}
and define $L = \{\by: \|\bw\| \le r\}$ and
$M=\{\by: \bw=\bA\bu+\bb \:\: \text{and} \:\: r = \bc^t\bu+d\}$,
then we can rephrase the problem as minimizing $\frac{1}{2}\|\bx-\bu\|^2$
subject to $\by \in L \cap M$.  This is a fairly typical set projection
problem except that the $\bw$ and $r$ components of $\by$ are missing 
in the loss function. 

Taking a cue from Example \ref{ex:lp}, we incorporate the affine 
constraints in the domain of the objective function. If we represent 
projection onto $L$ by
\begin{eqnarray*}
P \begin{pmatrix} \bw_k \\ r_k \end{pmatrix} & = & \begin{pmatrix} \tilde{\bw}_k \\ \tilde{r}_k \end{pmatrix},
\end{eqnarray*}
then the Lagrangian generated by the proximal distance algorithm amounts to
\begin{eqnarray*}
{\cal L} & = & \frac{1}{2}\|\bx-\bu\|^2+
\frac{\rho}{2}\left\|\begin{pmatrix} \bw -\tilde{\bw}_k \\ r - \tilde{r}_k \end{pmatrix}\right\|^2
+\blambda^t(\bA\bu+\bb-\bw)+\theta(\bc^t\bu+d-r).
\end{eqnarray*}
This gives rise to a system of three stationarity equations
\begin{eqnarray}
{\bf 0} & = & \bu-\bx + \bA^t \blambda+\theta \bc \label{lorentz1}\\
{\bf 0} & = & \rho(\bw - \tilde{\bw}_k)-\blambda \label{lorentz2} \\
0 & = & \rho(r - \tilde{r}_k) -\theta. \label{lorentz3}
\end{eqnarray}
Solving for the multipliers $\blambda$ and $\theta$ in equations
(\ref{lorentz2}) and (\ref{lorentz3}) and substituting their
values in equation (\ref{lorentz1}) yield
\begin{eqnarray*}
{\bf 0} & = & \bu-\bx +  \rho \bA^t (\bw-\tilde{\bw}_k) + \rho(r-\tilde{r}_k) \bc \\
& = & \bu-\bx + \rho \bA^t (\bA\bu+\bb-\tilde{\bw}_k) +\rho (\bc^t\bu+d-\tilde{r}_k) \bc.
\end{eqnarray*}
This leads to the MM update
\begin{eqnarray}
\label{eq:socp_update}
\bu_{k+1} & = & (\rho^{-1}\bI+ \bA^t\bA+ \bc \bc^t)^{-1}[\rho^{-1}\bx+\bA^t(\tilde{\bw}_k-\bb) 
+ (\tilde{r}_k -d)\bc].
\end{eqnarray}
The updates $\bw_{k+1} = \bA\bu_{k+1}+\bb$ and $r_{k+1} = \bc^t\bu_{k+1}+d$ 
follow from the constraints. 

Table \ref{tab:socp} compares the proximal distance algorithm to SCS and Gurobi. Echoing previous examples, 
we tailor the update schedule for $\rho$ differently for dense and sparse problems. Dense problems converge quickly and accurately when we set $\rho_0 = 1$ and double $\rho$ every 100 iterations. Sparse problems require a greater range and faster updates of $\rho$, so we set 
$\rho_0 = 0.01$ and then multiply $\rho$ by 2.5 every 10 iterations. For dense problems, it is clearly advantageous to cache the spectral decomposition of $\bA^t \bA + \bc \bc^t$ as suggested in Example \ref{ex:nqp}. In this regime, the proximal distance algorithm is as accurate as Gurobi and nearly as fast. SCS is comparable to Gurobi in speed but notably less accurate. 

With a large sparse constraint matrix $\bA$, extraction of
its spectral decomposition becomes prohibitive. If we let
$\bE = (\rho^{-1/2} \bI \:\: \bA^t \:\: \bc )$, then we must solve
a linear system of equations defined by the Gramian matrix $\bG = \bE \bE^t$.
There are three reasonable options for solving this system. 
The first relies on computing and caching a sparse Cholesky decomposition of $\bG$. 
The second computes the QR decomposition of the sparse matrix $\bE$.
The R part of the QR decomposition coincides with the Cholesky factor. Unfortunately, every time $\rho$ changes, the Cholesky or QR decomposition must be redone. 
The third option is the conjugate gradient algorithm. 
In our experience the QR decomposition offers superior stability and accuracy.
When $\bE$ is very sparse, the QR decomposition is often much faster than the Cholesky decomposition because it avoids forming the dense matrix $\bA^t \bA$. 
Even when only 5\% of the entries of $\bA$ are nonzero, 90\% of the entries of $\bA^t \bA$ can be nonzero. If exquisite accuracy is not a concern, then the conjugate gradient method provides the fastest update. Table 4 reflects this choice.

\begin{table}
    \centering
    \begin{tabular}{cccccccc}
        \toprule
        \multicolumn{2}{c}{Dimensions} & \multicolumn{3}{c}{Optima} & \multicolumn{3}{c}{CPU Seconds} \\
        \cmidrule(r){1-2} \cmidrule(r){3-5} \cmidrule(r){6-8}
        $m$ & $n$ & PD & SCS & Gurobi & PD & SCS & Gurobi \\
        \hline
        2    &     4 & 0.10598 & 0.10607 & 0.10598 & 0.0043 & 0.0103 & 0.0026 \\
        4    &     8 & 0.00000 & 0.00000 & 0.00000 & 0.0003 & 0.0009 & 0.0022 \\
        8    &    16 & 0.88988 & 0.88991 & 0.88988 & 0.0557 & 0.0011 & 0.0027 \\
        16   &    32 & 2.16514 & 2.16520 & 2.16514 & 0.0725 & 0.0012 & 0.0040 \\
        32   &    64 & 3.03855 & 3.03864 & 3.03853 & 0.0952 & 0.0019 & 0.0094 \\
        64   &   128 & 4.86894 & 4.86962 & 4.86895 & 0.1225 & 0.0065 & 0.0403 \\
        128  &   256 & 10.5863 & 10.5843 & 10.5863 & 0.1975 & 0.0810 & 0.0868 \\
        256  &   512 & 31.1039 & 31.0965 & 31.1039 & 0.5463 & 0.3995 & 0.3405 \\
        512  &  1024 & 27.0483 & 27.0475 & 27.0483 & 3.7667 & 1.6692 & 2.0189 \\
        1024 &  2048 & 1.45578 & 1.45569 & 1.45569 & 0.5352 & 0.3691 & 1.5489 \\
        2048 &  4096 & 2.22936 & 2.22930 & 2.22921 & 1.0845 & 2.4531 & 5.5521 \\
        4096 &  8192 & 1.72306 & 1.72202 & 1.72209 & 3.1404 & 17.272 & 15.204 \\
        8192 & 16384 & 5.36191 & 5.36116 & 5.36144 & 13.979 & 133.25 & 88.024 \\
        \bottomrule
    \end{tabular}
    \caption{CPU times and optima for the second-order cone projection. Here $m$ is the number of constraints, $n$ is the number of variables, PD is the accelerated proximal distance algorithm, SCS is the Splitting Cone Solver, and Gurobi is the Gurobi solver. After $m = 512$ the constraint matrix $\boldsymbol{A}$ is initialized with sparsity level 0.01.}
    \label{tab:socp}
\end{table}

\subsection{Copositive Matrices}
\label{ex:copos}

A symmetric matrix $\bM$ is copositive if
its associated quadratic form $\bx^t \bM \bx$ is nonnegative
for all $\bx \ge {\bf 0}$. Copositive matrices
find applications in numerous branches of the mathematical
sciences \citep{berman1994nonnegative}.  All positive
semidefinite matrices and all matrices with nonnegative
entries are copositive. The variational index
\begin{eqnarray*}
\mu(\bM) & = & \min_{\|\bx\| = 1, \: \bx \ge {\bf 0}} \bx^t \bM \bx
\end{eqnarray*}
is one key to understanding copositive matrices \citep{hiriart2010variational}. 
The constraint set $S$ is the intersection of the unit sphere and the 
nonnegative cone $\mathbb{R}_{+}^{n}$. Projection of an external point $\by$ onto 
$S$ splits into three cases. When all components of $\by$ are negative, then
$P_S(\by) = \be_i$, where $y_i$ is the least negative component 
of $\by$, and $\be_i$ is the standard unit vector along coordinate
direction $i$. The origin $\bf 0$ is equidistant from all points of $S$.
If any component of $\by$ is positive, then the projection is constructed
by setting the negative components of $\by$ equal
to 0, and standardizing the truncated version of $\by$ to have Euclidean norm 1.

As a test case for the proximal distance algorithm, consider the Horn matrix 
\citep{hall1963copositive}
\begin{eqnarray*}
\bM & = & \left[ \!\!\begin{array}{rrrrrr}
1 & -1 & 1 & 1 & -1 \\
-1 & 1 & -1 & 1 & 1 \\
1 & -1 & 1 & -1 & 1 \\
1 & 1 & -1 & 1 & -1 \\
-1 & 1 & 1 & -1 & 1 \end{array} \! \right].
\end{eqnarray*}
The value $\mu(\bM) = 0$ is attained for the vectors
$\frac{1}{\sqrt{2}}(1,1,0,0,0)^t$, $\frac{1}{\sqrt{6}}(1,2,1,0,0)^t$,
and equivalent vectors with their entries permuted. Matrices in higher 
dimensions with the same Horn pattern of 1's and -1's are 
copositive as well \citep{johnson2008constructing}. A Horn matrix of odd dimension
cannot be written as a positive semidefinite matrix, a nonnegative matrix,
or a sum of two such matrices. 

The proximal distance algorithm minimizes the criterion
\begin{eqnarray*}
g(\bx \mid \bx_k) & = & \frac{1}{2} \bx^t \bM \bx+\frac{\rho}{2}\|\bx-P_S(\bx_k)\|^2 
\end{eqnarray*}
and generates the updates
\begin{eqnarray*}
\bx_{k+1} & = & (\bM+\rho \bI)^{-1} \rho P_S(\bx_k).
\end{eqnarray*}
It takes a gentle tuning schedule to get decent results. The choice 
$\rho_k = 1.2^k$ converges in 600 to 700
iterations from random starting points and reliably yields 
objective values below $10^{-5}$ for Horn matrices.  The computational burden
per iteration is significantly eased by exploiting the
cached spectral decomposition of $\bM$. Table \ref{tab:copos-horn}
compares the performance of the proximal distance algorithm to the
Mosek solver on a range of Horn matrices. Mosek uses semidefinite programming to decide whether
$\bM$ can be decomposed into a sum of a positive semidefinite matrix 
and a nonnegative matrix. If not, Mosek declares the problem infeasible.
Nesterov acceleration improves the final loss for the proximal distance
algorithm, but it does not decrease overall computing time.

\begin{table}
	\centering
	\begin{tabular}{crrrrrr}
        \toprule
		\multicolumn{1}{c}{Dimension} & \multicolumn{3}{c}{Optima} & 
		\multicolumn{3}{c}{CPU Seconds} \\ 
		\cmidrule(r){1-1} \cmidrule(r){2-4} \cmidrule(r){5-7}
		$n$ & \multicolumn{1}{c}{PD} & \multicolumn{1}{c}{aPD} & 
		\multicolumn{1}{c}{Mosek} & \multicolumn{1}{c}{PD} & 
		\multicolumn{1}{c}{aPD} & \multicolumn{1}{c}{Mosek}\\
        \hline
		4 & 0.000000 & 0.000000 & feasible & 0.5555 & 0.0124 & 2.7744\\
		5 & 0.000000 & 0.000000 & infeasible & 0.0039 & 0.0086 & 0.0276\\
		8 & 0.000021 & 0.000000 & feasible & 0.0059 & 0.0083 & 0.0050\\
		9 & 0.000045 & 0.000000 & infeasible & 0.0055 & 0.0072 & 0.0082\\
		16 & 0.000377 & 0.000001 & feasible & 0.0204 & 0.0237 & 0.0185\\
		17 & 0.000441 & 0.000001 & infeasible & 0.0204 & 0.0378 & 0.0175\\
		32 & 0.001610 & 0.000007 & feasible & 0.0288 & 0.0288 & 0.1211\\
		33 & 0.002357 & 0.000009 & infeasible & 0.0242 & 0.0346 & 0.1294\\
		64 & 0.054195 & 0.000026 & feasible & 0.0415 & 0.0494 & 3.6284\\
		65 & 0.006985 & 0.000026 & infeasible & 0.0431 & 0.0551 & 2.7862\\
		\bottomrule
	\end{tabular}
	\caption{CPU times (seconds) and optima for approximating the Horn variational
	index of a Horn matrix. Here $n$ is the size of Horn matrix, PD is the proximal distance algorithm, aPD is the accelerated proximal distance algorithm, and Mosek is the Mosek solver.}
	\label{tab:copos-horn}
\end{table}

Testing for copositivity is challenging because neither the loss function 
nor the constraint set is convex. The proximal distance
algorithm offers a fast screening device for checking 
whether a matrix is copositive. On random $1000 \times 1000$ 
symmetric matrices $\bM$, the method invariably returns a negative 
index in less than two seconds of computing time.  Because the vast majority 
of symmetric matrices are not copositive, accurate estimation of the minimum 
is not required. Table \ref{tab:copos-sym} summarizes a few random trials
with lower-dimensional symmetric matrices. In higher dimensions, Mosek becomes 
non-competitive, and Nesterov acceleration is of dubious value.  

\begin{table}
	\centering
	\begin{tabular}{crrrrrr}
        \toprule
		\multicolumn{1}{c}{Dimension} & \multicolumn{3}{c}{Optima} & 
		\multicolumn{3}{c}{CPU Seconds} \\ 
		\cmidrule(r){1-1}  \cmidrule(r){2-4} \cmidrule(r){5-7}
		$n$ & \multicolumn{1}{c}{PD} & \multicolumn{1}{c}{aPD} & 
		\multicolumn{1}{c}{Mosek} & \multicolumn{1}{c}{PD} & 
		\multicolumn{1}{c}{aPD} & Mosek\\
        \hline
		4 & -0.391552 & -0.391561 & infeasible & 0.0029 & 0.0031 & 0.0024\\
		8 & -0.911140 & -2.050316 & infeasible & 0.0037 & 0.0044 & 0.0045\\
		16 & -1.680697 & -1.680930 & infeasible & 0.0199 & 0.0272 & 0.0062\\
		32 & -2.334520 & -2.510781 & infeasible & 0.0261 & 0.0242 & 0.0441\\
		64 & -3.821927 & -3.628060 & infeasible & 0.0393 & 0.0437 & 0.6559\\
		128 & -5.473609 & -5.475879 & infeasible & 0.0792 & 0.0798 & 38.3919\\
		256 & -7.956365 & -7.551814 & infeasible & 0.1632 & 0.1797 & 456.1500\\
		\bottomrule
	\end{tabular}
	\caption{CPU times and optima for testing the copositivity of random symmetric matrices. Here $n$ is the size of matrix, PD is the proximal distance algorithm, 
aPD is the accelerated proximal distance algorithm, and Mosek is the Mosek solver.}
	\label{tab:copos-sym}
\end{table}

\subsection{Linear Complementarity Problem}
\label{ex:lcp}

The linear complementarity problem \citep{murty1988linear} consists of finding vectors
$\bx$ and $\by$ with nonnegative components such that $\bx^t\by = 0$
and $\by = \bA \bx+\bb$ for a given square matrix $\bA$ and vector
$\bb$. The natural loss function is $\frac{1}{2}\|\by-\bA\bx-\bb\|^2$.
To project a vector pair $(\bu,\bv)$ onto the nonconvex constraint
set, one considers each component pair $(u_i,v_i)$
in turn. If $u_i \ge \max\{v_i,0\}$, then the nearest pair $(\bx,\by)$ has
components $(x_i,y_i) = (u_i,0)$. If $v_i \ge \max\{u_i,0\}$, then 
the nearest pair has components $(x_i,y_i) = (0,v_i)$. Otherwise,
$(x_i,y_i) = (0,0)$.  At each iteration the proximal
distance algorithm minimizes the criterion
\begin{eqnarray*}
\frac{1}{2}\|\by-\bA\bx-\bb\|^2+ \frac{\rho}{2}\|\bx - \tilde{\bx}_k\|^2
+\frac{\rho}{2}\|\by-\tilde{\by}_k\|^2,
\end{eqnarray*}
where $(\tilde{\bx}_k,\tilde{\by}_k)$ is the projection of $(\bx_k,\by_k)$
onto the constraint set. The stationarity equations become
\begin{eqnarray*}
{\bf 0} & = & -\bA^t(\by-\bA\bx-\bb) +\rho(\bx-\tilde{\bx}_k) \\
{\bf 0} & = & \by-\bA\bx-\bb+\rho(\by-\tilde{\by}_k) .
\end{eqnarray*}
Substituting the value of $\by$ from the second equation into
the first equation leads to the updates
\begin{eqnarray}
\bx_{k+1} & = & [(1+\rho)\bI+\bA^t\bA]^{-1}
[\bA^t(\tilde{\by}_k-\bb)+(1+\rho)\tilde{\bx}_k] \label{lc_system} \\
\by_{k+1} & = & \frac{1}{1+\rho}(\bA\bx_{k+1}+\bb)+\frac{\rho}{1+\rho}\tilde{\by}_k.
\nonumber
\end{eqnarray}
The linear system (\ref{lc_system}) can be solved in low to moderate dimensions
by computing and caching the spectral decomposition of $\bA^t \bA$ and in
high dimensions by the conjugate gradient method. Table \ref{tab:lincomp} compares
the performance of the proximal gradient algorithm to the Gurobi solver
on some randomly generated problems.

\begin{table}
	\centering
	\begin{tabular}{crrrr}
        \toprule
		\multicolumn{1}{c}{Dimension} & \multicolumn{2}{c}{Optima} & 
		\multicolumn{2}{c}{CPU Seconds} \\ 
		\cmidrule(r){1-1} \cmidrule(r){2-3} \cmidrule(r){4-5}
		$n$ & \multicolumn{1}{c}{PD} & \multicolumn{1}{c}{Mosek} & 
		\multicolumn{1}{c}{PD} & Mosek \\
        \hline
		4 & 0.000000 & 0.000000 & 0.0230 & 0.0266\\
		8 & 0.000000 & 0.000000 & 0.0062 & 0.0079\\
		16 & 0.000000 & 0.000000 & 0.0269 & 0.0052\\
		32 & 0.000000 & 0.000000 & 0.0996 & 0.4303\\
		64 & 0.000074 & 0.000000 & 2.6846 & 360.5183\\
		\bottomrule
	\end{tabular}
	\caption{CPU times (seconds) and optima for the linear complementarity problem with randomly generated data. Here $n$ is the size of matrix, PD is the accelerated proximal distance algorithm, and Gurobi is the Gurobi solver.}
	\label{tab:lincomp}
\end{table}

\subsection{Sparse Principal Components Analysis}
\label{ex:spca}

Let $\bX$ be an $n \times p$ data matrix gathered on $n$ cases and $p$ predictors.
Assume the columns of $\bX$ are centered to have mean 0. Principal component analysis
(PCA) \citep{hotelling1933,pearson1901} operates on the sample covariance matrix $\bS= \frac{1}{n}\bX^t \bX$. 
Here we formulate a proximal distance algorithm for sparse PCA (SPCA),
which has attracted substantial interest in the machine learning community \citep{berthet2013optimal,berthet2013complexity,daspremont2007,johnstonelu2009,journee2010,wittentibshiranihastie2009,zouhastietibshirani2006}.
According to a result of Ky Fan \citep{fan1949theorem},
the first $q$ principal components (PCs) $\bu_1,\ldots,\bu_q$ can be extracted 
by maximizing the function $\tr(\bU^t\bS \bU)$ subject to the matrix
constraint $\bU^t \bU = \bI_q$, where $\bu_i$ is the $i$th column of
the $p \times q$ matrix $\bU$. This constraint set is called a Stiefel manifold.
One can impose sparsity by insisting that any given column $\bu_i$ have 
at most $r$ nonzero entries. Alternatively, one can require the entire matrix $\bU$ to have at most $r$ nonzero entries. The latter choice permits sparsity to be distributed non-uniformly
across columns.  

Extraction of sparse PCs is difficult for three reasons. First,
the Stiefel manifold $M_q$ and both sparsity sets are nonconvex. Second,
the objective function is concave rather than convex. Third, 
there is no simple formula or algorithm for projecting onto the 
intersection of the two constraint sets. Fortunately, it is straightforward 
to project onto each separately. Let $P_{M_q}(\bU)$ denote 
the projection of $\bU$ onto the Stiefel manifold.
It is well known that $P_{M_q}(\bU)$ can be calculated by extracting
a partial singular value decomposition $\bU = \bV \bSigma \bW^t$ of $\bU$
and setting $P_{M_q}(\bU) = \bV \bW^t$ \citep{golub2012matrix}. 
Here $\bV$ and $\bW$ are orthogonal matrices of dimension $p \times q$ and $q \times q$, respectively, and $\bSigma$ is a diagonal matrix of dimension $q \times q$. 
Let $P_{S_r}(\bU)$ denote the projection of $\bU$ onto the sparsity set
\begin{eqnarray*}
S_r & = & \{\bV: v_{ij} \ne 0 \; \text{for at most $r$ entries of each column $\bv_i$}\}.
\end{eqnarray*}
Because $P_{S_r}(\bU)$ operates column by column, it suffices
to project each column vector $\bu_i$ to sparsity. This entails nothing more
than sorting the entries of $\bu_i$ by magnitude, saving the $r$ largest,
and sending the remaining $p-r$ entries to 0. If the entire matrix $\bU$ must
have at most $r$ nonzero entries, then $\bU$ can be treated as a concatenated vector
during projection.

The key to a good algorithm is to incorporate the Stiefel 
constraints into the domain of the objective function \citep{kiers1990majorization,kiers1992minimization}
and the sparsity constraints into the distance penalty. Thus, 
we propose decreasing the criterion 
\begin{eqnarray*}
f(\bU) & = & -\frac{1}{2}\tr(\bU^t\bS \bU)+\frac{\rho}{2}\dist(\bU,S_r)^2 .
\end{eqnarray*}
at each iteration subject to the Stiefel constraints. The loss can be majorized via
\begin{eqnarray*}
 -\frac{1}{2}\tr(\bU^t\bS \bU) & = & 
-\frac{1}{2}\tr[(\bU-\bU_k)^t\bS (\bU-\bU_k)]
-\tr(\bU^t \bS \bU_k) + \frac{1}{2}\tr(\bU_k^t \bS \bU_k) \\
& \le & -\tr(\bU^t \bS \bU_k) + \frac{1}{2}\tr(\bU_k^t \bS \bU_k) 
\end{eqnarray*}
because $\bS$ is positive semidefinite. The penalty is majorized by
\begin{eqnarray*}
\frac{\rho}{2}\dist(\bU,S_r)^2 & \le & - \rho \tr[\bU^tP_{S_r}(\bU_k)]+c_k
\end{eqnarray*}
up to an irrelevant constant $c_k$ since the squared Frobenius norm 
satisfies the relation $\|\bU^t \bU\|_F^2 = q$ on the Stiefel manifold. It now follows that $f(\bU)$ 
is majorized by
\begin{eqnarray*}
\frac{1}{2}\|\bU- \bS \bU_k - \rho P_{S_r}(\bU_k)\|_F^2 
\end{eqnarray*}
up to an irrelevant constant. Accordingly, the Stiefel projection 
\begin{eqnarray*}
\bU_{k+1} & = & P_{M_q}[\bS \bU_k + \rho P_{S_r}(\bU_k)]
\end{eqnarray*}
provides the next MM iterate. 

\begin{figure}
    \includegraphics[width=0.9\textwidth]{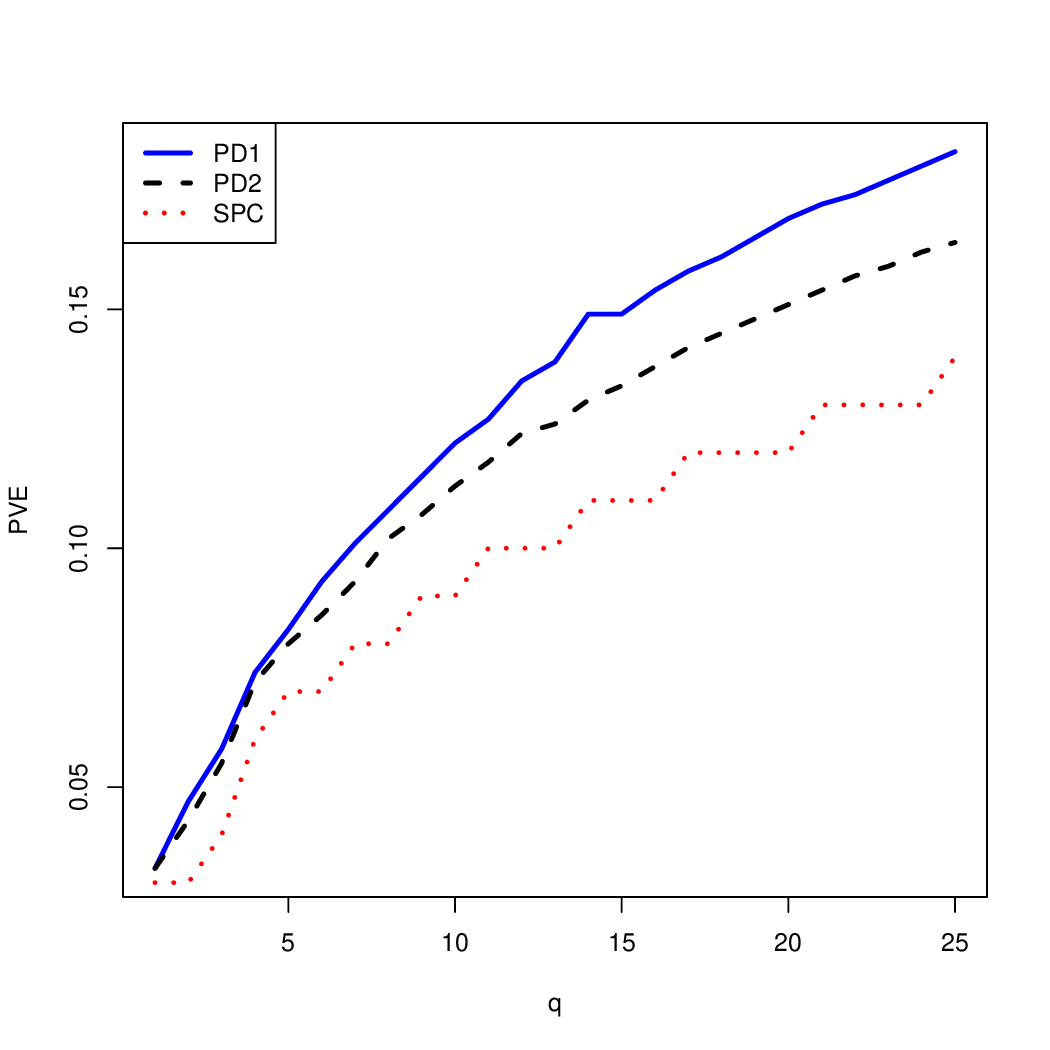}
    \caption{Proportion of variance explained by $q$ PCs for each algorithm. Here PD1 is the accelerated proximal distance algorithm enforcing matrix sparsity, PD2 is the accelerated proximal distance algorithm enforcing column-wise sparsity, and SPC is the orthogonal sparse PCA method from \texttt{PMA}.} 
    \label{fig:pve}
\end{figure}

\begin{figure}
    \includegraphics[width=0.9\textwidth]{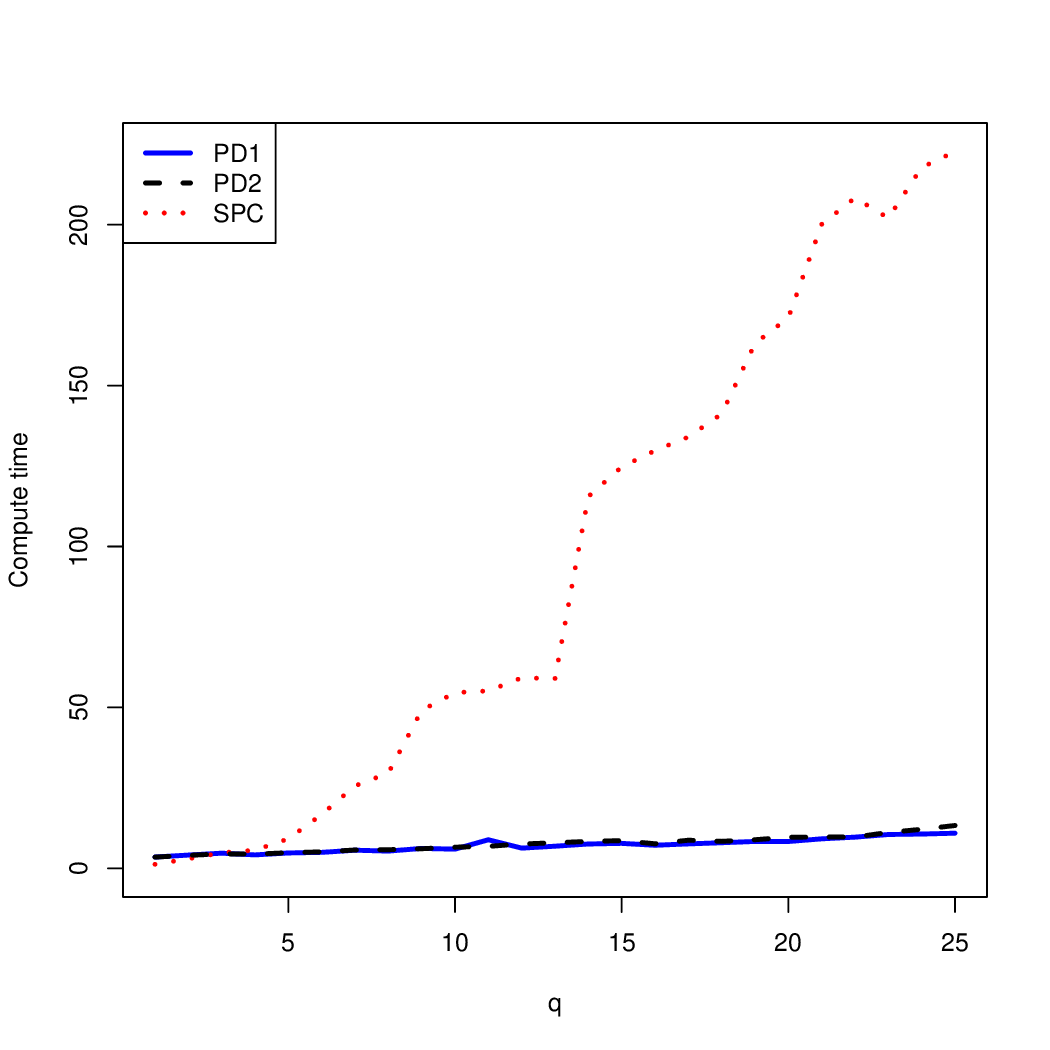}
    \caption{Computation times for $q$ PCs for each algorithm. Here PD1 is the accelerated proximal distance algorithm enforcing matrix sparsity, PD2 is the accelerated proximal distance algorithm  enforcing column-wise sparsity, and SPC is the orthogonal sparse PCA method from \texttt{PMA}.} 
    \label{fig:times}
\end{figure}

Figures \ref{fig:pve} and \ref{fig:times} compare the proximal distance algorithm to the SPC function from the R package \texttt{PMA} \citep{wittentibshiranihastie2009}.
The breast cancer  data from \texttt{PMA} provide the data matrix $\bX$.
The data consist of $p = 19,672$ RNA measurements on $n = 89$ patients. The two
figures show computation times and the proportion of variance explained (PVE) by the 
$p \times q$ loading matrix $\bU$. For sparse PCA, PVE is defined as 
$\tr(\bX_{q}^t \bX_{q}) / \tr(\bX^{t} \bX)$, where 
$\bX_{q} = \bX \bU(\bU^t \bU)^{-1} \bU^t$ \citep{shenhuang2008}.
When the loading vectors of $\bU$ are orthogonal, this criterion reduces to the familiar definition $\tr(\bU^t \bX^t \bX \bU) / \tr(\bX^t \bX)$ of PVE for ordinary PCA.
The proximal distance algorithm enforces either matrix-wise or column-wise
sparsity. In contrast, SPC enforces only column-wise sparsity via the constraint
$\|\bu_i\|_1 \le c$ for each column $\bu_i$ of $\bU$. We take $c=8$. The number of nonzeroes 
per loading vector output by SPC dictates the sparsity level for the column-wise
version of the proximal distance algorithm. Summing these counts across all columns
dictates the sparsity level for the matrix version of the proximal distance
algorithm. 

Figures \ref{fig:pve} and \ref{fig:times} demonstrate the superior PVE and computational speed of both proximal distance algorithms versus SPC. 
The type of projection does not appear to affect the computational performance of the proximal distance algorithm, as both versions scale equally well with $q$. 
However, the matrix projection, which permits the algorithm to more freely assign nonzeroes to the loadings, attains better PVE than the more restrictive column-wise projection.
For both variants of the proximal distance algorithm, Nesterov acceleration improves both fitting accuracy and computational speed, especially as the number of PCs $q$ increases.

\section{Discussion}
\label{sec:discussion}

The proximal distance algorithm applies to a host of problems. 
In addition to the linear and quadratic programming examples considered here, our previous paper \citep{lange2014proximal} derives and tests algorithms for binary piecewise-linear programming, $\ell_0$ regression, matrix completion \citep{cai10,candes10,chen2012matrix,mazumder2010spectral}, and sparse precision matrix estimation \citep{friedman08}. 
Other potential applications immediately come to mind. An integer linear program in standard form can be expressed as minimizing $\bc^t \bx$  subject to
$\bA \bx + \bs = \bb$, $\bs \ge {\bf 0}$, and $\bx \in \mathbb{Z}^p$. The latter
two constraints can be combined in a single constraint for which projection is trivial.
The affine constraints should be folded into the domain of the objective. Integer
programming is NP hard, so that the proximal distance algorithm just sketched is merely heuristic. Integer linear programming includes traditional NP hard problems such as the traveling salesman problem, the vertex cover problem, set packing, and Boolean satisfiability. It will be interesting to see if the proximal distance principle is competitive in meeting these challenges. Our experience with the closest lattice point problem \citep{agrell2002closest} and the eight queens problem suggests that the proximal distance algorithm can be too greedy for hard combinatorial problems.  The nonconvex problems solved in this paper are in some vague sense easy 
combinatorial problems.

The behavior of a proximal distance algorithm depends critically on a sensible tuning schedule for increasing $\rho$.
Starting $\rho$ too high puts too much stress on satisfying the constraints. 
Incrementing $\rho$ too quickly causes the algorithm to veer off the solution path guaranteed by the penalty method.  
Given the chance of roundoff error even with double precision arithmetic, it is unwise to take $\rho$ all the way to $\infty$. 
Trial and error can help in deciding whether a given class of problems will benefit from an aggressive update schedule
and strict or loose convergence criteria.
In problems with little curvature such as linear programming, more conservative updates are probably prudent. 
The linear programming, closest kinship matrix, and SPCA problems document the value of folding constraints into the domain of the loss. 
In the same spirit it is wise to minimize the number of constraints. A single penalty for projecting onto the intersection of two constraint sets is almost always preferable to the sum of two penalties for their separate projections. Exceptions to this rule obviously occur when projection onto the intersection is intractable. The integer linear programming problem mentioned previously illustrates these ideas. 

Our earlier proximal distance algorithms ignored acceleration. In
many cases the solutions produced had very low accuracy. The realization that
convex proximal distance algorithms can be phrased as proximal gradient
algorithms convinced us to try Nesterov acceleration. We now do this routinely
on the subproblems with $\rho$ fixed. This typically forces tighter path
following and a reduction in overall computing times. Our examples generally 
bear out the contention that Nesterov acceleration is useful in nonconvex 
problems \citep{Ghadimi2015}. It is noteworthy that the value of acceleration often lies in improving the quality of a solution as much as it does in increasing the rate of convergence. Of course, acceleration cannot prevent convergence to an inferior local minimum. 

On both convex and nonconvex problems, proximal distance algorithms enjoy global convergence guarantees. On nonconvex problems, one must confine attention to subanalytic sets and subanalytic functions. This minor restriction is not a handicap in practice.  Determining local convergence rates is a more vexing issue. For convex problems, we review existing theory buttressing an $O(\rho k^{-1})$ sublinear rate. Better results require restrictive smoothness assumptions on both the objective function and the constraint sets. When $f(\bx)$ is convex, but the constraint sets are nonconvex, proximal distance algorithms reduce to concave-convex programming. \cite{le2009convergence} attack convergence rates for concave-convex programs.

We hope readers will sense the potential of the proximal distance principle. This simple idea offers insight into many existing algorithms and a straightforward path in devising new ones.  Effective proximal and projection operators usually spell the difference between success and failure. The number and variety of such operators is expanding quickly as the field of optimization relinquishes it fixation on convexity.  The current paper research leaves many open questions about tuning schedules, rates of convergence, and acceleration in the face of nonconvexity. We welcome the contributions of other mathematical scientists in unraveling these mysteries and in inventing new proximal distance algorithms.

\section*{\center Acknowledgments}
We thank Joong-Ho Won for many insightful discussions. In partiular, he pointed out the utility of the least squares criterion (\ref{johann_credit}).
Hua Zhou and Kenneth Lange were supported by grants from the National Human Genome Research Institute (HG006139) and the National Institute of General Medical Sciences (GM053275).
Kevin Keys was supported by a National Science Foundation Graduate Research Fellowship (DGE-0707424), a Predoctoral Training Grant (HG002536) from the National Human Genome Research Institute, and by a National Heart, Lung, Blood Institute grant (R01HL135156).

\bibliographystyle{plain}
\bibliography{ProximalDistance}

\begin{thebibliography}{10}

\bibitem{agrell2002closest}
Erik Agrell, Thomas Eriksson, Alexander Vardy, and Kenneth Zeger.
\newblock Closest point search in lattices.
\newblock {\em IEEE Transactions on Information Theory}, 48(8):2201--2214,
  2002.

\bibitem{alizadeh2003second}
Farid Alizadeh and Donald Goldfarb.
\newblock Second-order cone programming.
\newblock {\em Mathematical {P}rogramming}, 95:3--51, 2003.

\bibitem{attouch2010proximal}
H{\'e}dy Attouch, J{\'e}r{\^o}me Bolte, Patrick Redont, and Antoine Soubeyran.
\newblock Proximal alternating minimization and projection methods for
  nonconvex problems: An approach based on the kurdyka-{\l}ojasiewicz
  inequality.
\newblock {\em Mathematics of Operations Research}, 35(2):438--457, 2010.

\bibitem{bauschke2011convex}
Heinz~H Bauschke and Patrick~L Combettes.
\newblock {\em Convex {A}nalysis and {M}onotone {O}perator {T}heory in
  {H}ilbert {S}paces}.
\newblock Springer, 2011.

\bibitem{beck2009fast}
Amir Beck and Marc Teboulle.
\newblock A fast iterative shrinkage-thresholding algorithm for linear inverse
  problems.
\newblock {\em SIAM Journal on Imaging Sciences}, 2(1):183--202, 2009.

\bibitem{beltrami1970algorithmic}
Edward~J Beltrami.
\newblock {\em An Algorithmic Approach to Nonlinear Analysis and Optimization}.
\newblock Academic Press, 1970.

\bibitem{bemporad2018semidefiniteQP}
Alberto Bemporad.
\newblock A numerically stable solver for positive semidefinite quadratic
  programs based on nonnegative least squares.
\newblock {\em IEEE Transactions on Automatic Control}, 63(2):525--531, 2018.

\bibitem{berman1994nonnegative}
Abraham Berman and Robert~J Plemmons.
\newblock {\em Nonnegative Matrices in the Mathematical Sciences}.
\newblock Classics in Applied Mathematics. Society for Industrial and Applied
  Mathematics, 1994.

\bibitem{berthet2013complexity}
Quentin Berthet and Philippe Rigollet.
\newblock Complexity theoretic lower bounds for sparse principal component
  detection.
\newblock In {\em Conference on Learning Theory}, pages 1046--1066, 2013.

\bibitem{berthet2013optimal}
Quentin Berthet and Philippe Rigollet.
\newblock Optimal detection of sparse principal components in high dimension.
\newblock {\em The Annals of Statistics}, 41(4):1780--1815, 2013.

\bibitem{bierstone1988semianalytic}
Edward Bierstone and Pierre~D Milman.
\newblock Semianalytic and subanalytic sets.
\newblock {\em Publications Math{\'e}matiques de l'Institut des Hautes
  {\'E}tudes Scientifiques}, 67(1):5--42, 1988.

\bibitem{bochnak2013real}
Jacek Bochnak, Michel Coste, and Marie-Fran{\c{c}}oise Roy.
\newblock {\em Real algebraic geometry}, volume~36.
\newblock Springer Science \& Business Media, 2013.

\bibitem{bolte2007lojasiewicz}
J{\'e}r{\^o}me Bolte, Aris Daniilidis, and Adrian Lewis.
\newblock The {\l}ojasiewicz inequality for nonsmooth subanalytic functions
  with applications to subgradient dynamical systems.
\newblock {\em SIAM Journal on Optimization}, 17(4):1205--1223, 2007.

\bibitem{borwein06}
Jonathan~M Borwein and Adrian~S Lewis.
\newblock {\em Convex Analysis and Nonlinear Optimization: Theory and
  Examples}.
\newblock CMS Books in Mathematics. Springer, New York, 2nd edition, 2006.

\bibitem{boyd2009convex}
Stephen Boyd and Lieven Vandenberghe.
\newblock {\em Convex {O}ptimization}.
\newblock Cambridge {U}niversity {P}ress, 2009.

\bibitem{boyle1986method}
James~P Boyle and Richard~L Dykstra.
\newblock A method for finding projections onto the intersection of convex sets
  in {H}ilbert spaces.
\newblock In {\em {Advances in Order Restricted Statistical Inference}}, pages
  28--47. Springer, 1986.

\bibitem{bunea2010spades}
Florentina Bunea, Alexandre~B Tsybakov, Marten~H Wegkamp, and Adrian Barbu.
\newblock Spades and mixture models.
\newblock {\em The Annals of Statistics}, 38(4):2525--2558, 2010.

\bibitem{cai10}
Jian-Feng Cai, Emmanuel~J Cand\`{e}s, and Zuowei Shen.
\newblock A singular value thresholding algorithm for matrix completion.
\newblock {\em SIAM Journal on Optimization}, 20:1956--1982, 2010.

\bibitem{candes10}
Emmanuel~J Cand\`{e}s and Terence Tao.
\newblock The power of convex relaxation: near-optimal matrix completion.
\newblock {\em IEEE Transactions on Information Theory}, 56:2053--2080, 2010.

\bibitem{chen2012matrix}
Caihua Chen, Bingsheng He, and Xiaoming Yuan.
\newblock Matrix completion via an alternating direction method.
\newblock {\em IMA Journal of Numerical Analysis}, 32:227--245, 2012.

\bibitem{chi2014distance}
Eric~C Chi, Hua Zhou, and Kenneth Lange.
\newblock Distance majorization and its applications.
\newblock {\em Mathematical {P}rogramming {S}eries {A}}, 146:409--436, 2014.

\bibitem{combettes2011proximal}
Patrick~L Combettes and Jean-Christophe Pesquet.
\newblock Proximal splitting methods in signal processing.
\newblock In {\em {Fixed-Point Algorithms for Inverse Problems in Science and
  Engineering}}, pages 185--212. Springer, 2011.

\bibitem{courant1943variational}
Richard Courant.
\newblock Variational methods for the solution of problems of equilibrium and
  vibrations.
\newblock {\em Bulletin of the American Mathematical Society}, 49:1--23, 1943.

\bibitem{cui2018composite}
Ying Cui, Jong-Shi Pang, and Bodhisattva Sen.
\newblock Composite difference-max programs for modern statistical estimation
  problems.
\newblock {\em arXiv preprint arXiv:1803.00205}, 2018.

\bibitem{daspremont2007}
Alexandre D'Aspremont, Laurent~El Ghaoui, Michael~I Jordan, and Gert R~G
  Lanckriet.
\newblock A direct formulation for sparse {PCA} using semidefinite programming.
\newblock {\em SIAM Review}, 49(3):434--448, 2007.

\bibitem{DunningHuchetteLubin2015}
Iain Dunning, Joey Huchette, and Miles Lubin.
\newblock {JuMP}: {A} modeling language for mathematical optimization.
\newblock {\em arXiv:1508.01982 [math.OC]}, 2015.

\bibitem{fan1949theorem}
Ky~Fan.
\newblock On a theorem of {W}eyl concerning eigenvalues of linear
  transformations {I}.
\newblock {\em Proceedings of the National Academy of Sciences of the United
  States of America}, 35:652--655, 1949.

\bibitem{fong2011lsmr}
David Chin-Lung Fong and Michael Saunders.
\newblock Lsmr: An iterative algorithm for sparse least-squares problems.
\newblock {\em SIAM Journal on Scientific Computing}, 33(5):2950--2971, 2011.

\bibitem{friedman08}
Jerome Friedman, Trevor Hastie, and Robert Tibshirani.
\newblock {Sparse inverse covariance estimation with the graphical lasso}.
\newblock {\em Biostatistics}, 9:432--441, July 2008.

\bibitem{Ghadimi2015}
Saeed Ghadimi and Guanghui Lan.
\newblock Accelerated gradient methods for nonconvex nonlinear and stochastic
  programming.
\newblock {\em {M}athematical {P}rogramming}, 156(1):59--99, 2015.

\bibitem{golub2012matrix}
Gene~H Golub and Charles~F Van~Loan.
\newblock {\em Matrix Computations}.
\newblock JHU Press, 3 edition, 2012.

\bibitem{hall1963copositive}
Marshall Hall and Morris Newman.
\newblock Copositive and completely positive quadratic forms.
\newblock In {\em Mathematical Proceedings of the Cambridge Philosophical
  Society}, volume~59, pages 329--339. Cambridge Univ Press, 1963.

\bibitem{HeylenBurazerovicScheunders11SimplexLS}
Rob Heylen, D\v{z}evdet Burazerovi\'{c}, and Paul Scheunders.
\newblock Fully constrained least squares spectral unmixing by simplex
  projection.
\newblock {\em IEEE Transactions on Geoscience and Remote Sensing},
  49(11):4112--4122, Nov 2011.

\bibitem{higham2002computing}
Nicholas~J Higham.
\newblock Computing the nearest correlation matrix -— a problem from finance.
\newblock {\em IMA Journal of Numerical Analysis}, 22(3):329--343, 2002.

\bibitem{hiriart2010variational}
Jean-Baptiste Hiriart-Urruty and Alberto Seeger.
\newblock A variational approach to copositive matrices.
\newblock {\em SIAM {R}eview}, 52:593--629, 2010.

\bibitem{hotelling1933}
Harold Hotelling.
\newblock Analysis of a complex of statistical variables into principle
  components.
\newblock {\em Journal of Educational Psychology}, 24:417--441, 1933.

\bibitem{hunter04}
David~R. Hunter and Kenneth Lange.
\newblock A tutorial on {MM} algorithms.
\newblock {\em American Statistician}, 58:30--37, 2004.

\bibitem{johnson2008constructing}
Charles~R Johnson and Robert Reams.
\newblock Constructing copositive matrices from interior matrices.
\newblock {\em Electronic Journal of Linear Algebra}, 17:9--20, 2008.

\bibitem{johnstonelu2009}
Iain~M Johnstone and Arthur~Yu Lu.
\newblock On consistency and sparsity for principal components analysis in high
  dimensions.
\newblock {\em Journal of the American Statistical Association},
  104(486):682--693, 2009.

\bibitem{journee2010}
Michel Journ{\'e}e, Yurii Nesterov, Peter Richt{\'a}rik, and Rodolphe
  Sepulchre.
\newblock Generalized power method for sparse principal component analysis.
\newblock {\em Journal of Machine Learning Research}, 11:517--553, 2010.

\bibitem{kang2015global}
Yangyang Kang, Zhihua Zhang, and Wu-Jun Li.
\newblock On the global convergence of majorization minimization algorithms for
  nonconvex optimization problems.
\newblock {\em arXiv preprint arXiv:1504.07791}, 2015.

\bibitem{keshava2003survey}
Nirmal Keshava.
\newblock A survey of spectral unmixing algorithms.
\newblock {\em Lincoln laboratory journal}, 14(1):55--78, 2003.

\bibitem{kiers1990majorization}
Henk~AL Kiers.
\newblock Majorization as a tool for optimizing a class of matrix functions.
\newblock {\em Psychometrika}, 55:417--428, 1990.

\bibitem{kiers1992minimization}
Henk~AL Kiers and Jos~MF ten Berge.
\newblock Minimization of a class of matrix trace functions by means of refined
  majorization.
\newblock {\em Psychometrika}, 57:371--382, 1992.

\bibitem{krasnosel1955two}
Mark~Aleksandrovich Krasnosel'skii.
\newblock Two remarks on the method of successive approximations.
\newblock {\em Uspekhi Matematicheskikh Nauk}, 10(1):123--127, 1955.

\bibitem{kruger2003frechet}
A~Ya Kruger.
\newblock On fr{\'e}chet subdifferentials.
\newblock {\em Journal of Mathematical Sciences}, 116(3):3325--3358, 2003.

\bibitem{landweber1951iteration}
Louis Landweber.
\newblock An iteration formula for {F}redholm integral equations of the first
  kind.
\newblock {\em American {J}ournal of {M}athematics}, 73(3):615--624, 1951.

\bibitem{lange10}
Kenneth Lange.
\newblock {\em Optimization}.
\newblock Springer, New York, 2nd edition, 2010.

\bibitem{lange2016MM}
Kenneth Lange.
\newblock {\em {MM} {O}ptimization {A}lgorithms}.
\newblock SIAM, 2016.

\bibitem{lange2014proximal}
Kenneth Lange and Kevin~L Keys.
\newblock The proximal distance algorithm.
\newblock In {\em Proceedings of the 2014 International Congress of
  Mathematicians}, pages 95--116. Kyung Moon, August 2014.

\bibitem{le2009convergence}
HA~Le~Thi, VN~Huynh, and T~Pham~Dinh.
\newblock Convergence analysis of dc algorithm for dc programming with
  subanalytic data.
\newblock {\em Ann. Oper. Res. Technical Report, LMI, INSA-Rouen}, 2009.

\bibitem{lobo1998applications}
Miguel~Sousa Lobo, Lieven Vandenberghe, Stephen Boyd, and Herv{\'e} Lebret.
\newblock Applications of second-order cone programming.
\newblock {\em Linear Algebra and its Applications}, 284:193--228, 1998.

\bibitem{LubinDunningIJOC}
Miles Lubin and Iain Dunning.
\newblock Computing in operations research using {J}ulia.
\newblock {\em INFORMS Journal on Computing}, 27(2):238--248, 2015.

\bibitem{luenberger1984linear}
David~G Luenberger and Yinyu Ye.
\newblock {\em Linear and nonlinear programming}, volume~2.
\newblock Springer, 1984.

\bibitem{mairal2013optimization}
Julien Mairal.
\newblock Optimization with first-order surrogate functions.
\newblock In {\em International Conference on Machine Learning}, pages
  783--791, 2013.

\bibitem{mann1953mean}
W~Robert Mann.
\newblock Mean value methods in iteration.
\newblock {\em Proceedings of the American Mathematical Society},
  4(3):506--510, 1953.

\bibitem{markowitz1952portfolio}
Harry Markowitz.
\newblock Portfolio selection.
\newblock {\em The journal of finance}, 7(1):77--91, 1952.

\bibitem{mazumder2010spectral}
Rahul Mazumder, Trevor Hastie, and Robert Tibshirani.
\newblock Spectral regularization algorithms for learning large incomplete
  matrices.
\newblock {\em The Journal of Machine Learning Research}, 11:2287--2322, 2010.

\bibitem{moreau1962fonctions}
Jean-Jacques Moreau.
\newblock Fonctions convexes duales et points proximaux dans un espace
  hilbertien.
\newblock {\em Comptes Rendus de l'Acad{\'e}mie des Sciences de Paris A},
  255:2897--2899, 1962.

\bibitem{murty1988linear}
Katta~G Murty and Feng-Tien Yu.
\newblock {\em Linear {C}omplementarity, {L}inear and {N}onlinear
  {P}rogramming}.
\newblock Heldermann Verlag, West Berlin, 1988.

\bibitem{niculescu2006convex}
Constantin Niculescu and Lars-Erik Persson.
\newblock {\em Convex {F}unctions and their {A}pplications: a {C}ontemporary
  {A}pproach}.
\newblock Springer, 2006.

\bibitem{o2013conic}
Brendan O{'}Donoghue, Eric Chu, Neal Parikh, and Stephen Boyd.
\newblock Conic optimization via operator splitting and homogeneous self-dual
  embedding.
\newblock {\em Journel of Optimization Theory and Applications}, pages 1--27,
  2016.

\bibitem{paige1982algorithm}
Christopher~C Paige and Michael~A Saunders.
\newblock Algorithm 583: Lsqr: Sparse linear equations and least squares
  problems.
\newblock {\em ACM Transactions on Mathematical Software (TOMS)},
  8(2):195--209, 1982.

\bibitem{paigesaunders1982b}
Christopher~C Paige and Michael~A Saunders.
\newblock Algorithm 583: {LSQR}: Sparse linear equations and least squares
  problems.
\newblock {\em ACM Transactions on Mathematical Software (TOMS)},
  8(2):195--209, 1982.

\bibitem{paigesaunders1982a}
Christopher~C Paige and Michael~A Saunders.
\newblock {LSQR}: An algorithm for sparse linear equations and sparse least
  squares.
\newblock {\em ACM Transactions on Mathematical Software (TOMS)}, 8(1):43--71,
  1982.

\bibitem{parikhboyd2013}
Neal Parikh and Stephen Boyd.
\newblock Proximal algorithms.
\newblock {\em Foundations and Trends in Optimization}, 1(3):123--231, 2013.

\bibitem{pearson1901}
Karl Pearson.
\newblock On lines and planes of closest fit to systems of points in space.
\newblock {\em Philosophical Magazine}, 2(11):559--572, 1901.

\bibitem{shenhuang2008}
Haipeng Shen and Jianhua~Z. Huang.
\newblock Sparse principal component analysis via regularized low rank matrix
  approximation.
\newblock {\em Journal of Multivariate Analysis}, 99:1015--1034, 2008.

\bibitem{su2014differential}
Weijie Su, Stephen Boyd, and Emmanuel Cand{\`e}s.
\newblock A differential equation for modeling {N}esterov's accelerated
  gradient method: Theory and insights.
\newblock In {\em Advances in Neural Information Processing Systems}, pages
  2510--2518, 2014.

\bibitem{waechterbiegler2005}
Andreas W{\"a}chter and Lorenz~T Biegler.
\newblock Line search filter methods for nonlinear programming: Motivation and
  global convergence.
\newblock {\em SIAM Journal on Optimization}, 16(1):1--31, 2005.

\bibitem{waechterbiegler2006}
Andreas W{\"a}chter and Lorenz~T Biegler.
\newblock On the implementation of an interior-point filter line-search
  algorithm for large-scale nonlinear programming.
\newblock {\em Mathematical {P}rogramming}, 106(1):25--57, 2006.

\bibitem{wittentibshiranihastie2009}
Daniela~M Witten, Robert Tibshirani, and Trevor Hastie.
\newblock A penalized matrix decomposition, with applications to sparse
  principal components and canonical correlation analysis.
\newblock {\em Biostatistics}, 10(3):515--534, 2009.

\bibitem{xu2017generalized}
Jason Xu, Eric Chi, and Kenneth Lange.
\newblock Generalized linear model regression under distance-to-set penalties.
\newblock In {\em Advances in Neural Information Processing Systems}, pages
  1385--1394, 2017.

\bibitem{zouhastietibshirani2006}
Hui Zou, Trevor Hastie, and Robert Tibshirani.
\newblock Sparse principal components analysis.
\newblock {\em Journal of {C}omputational and {G}raphical {S}tatistics},
  15(2):262--282, 2006.

\end{thebibliography}
\end{document}